 \documentclass[a4,12pt]{birkmult}

  \usepackage{latexsym}
  \usepackage{amsmath}   %needed for \begin{align}... \end{align} environment
  \usepackage{amsfonts}
  \usepackage{amssymb}
  \usepackage{graphicx}
  \usepackage{amsthm}
  \usepackage[noadjust]{cite}

  \newcommand{\G}{C \kern -0.1em \ell}
  
  \newcommand{\R}{\mathbb{R}}

  \newcommand{\qi}{\ensuremath{\mbox{\boldmath $i$}}}

  \newcommand{\qj}{\ensuremath{\mbox{\boldmath $j$}}}
  
  \newcommand{\qk}{\ensuremath{\mbox{\boldmath $k$}}}
  
  \newcommand{\be}{\begin{equation}}
  \newcommand{\ee}{\end{equation}} % Does not work. Don't know why.
  \newcommand{\bv}[1]{\mathbf{#1}}
  \newcommand{\vect}[1]{\vec{#1}}

  \newcommand{\vb}{\vec{b}}
  \newcommand{\bvc}{\underline{c}} 
  \newcommand{\alphap}{{\alpha^{\prime}}} 
  \newcommand{\vbp}{\vec{b}^{\prime}} 
  \newcommand{\vbpsq}{\vec{b}^{\prime 2}} 
  \newcommand{\bvcp}{{{\underline{c}}^{\prime}}} 
  \newcommand{\betap}{{\beta^{\prime}}}

  \newtheorem{example}{Example}

  \newtheorem{thm}{Theorem}[section]

  \theoremstyle{definition}
  \newtheorem{defn}[thm]{Definition}
  \theoremstyle{remark}

  \newcommand{\topline}{\hline\noalign{\smallskip}}
  \newcommand{\midline}{\noalign{\smallskip}\hline\noalign{\smallskip}}
  \newcommand{\bottomline}{\noalign{\smallskip}\hline\noalign{\smallskip}}

\newcommand{\ed}{\end{document}} 
\newcommand{\e}[1]{\Vec{e}_#1}
 
\begin{document}

%%% Old title: \title{Geometric roots of $-1$}

\title[Geometric Roots of $-1$ in Clifford Algebras $\G_{p,q}$ with $p+q \leq 4$]
{Geometric Roots of $-1$ in Clifford
%(Geometric)
Algebras\\ $\G_{p,q}$ with $p+q \leq 4$}

\author[E. Hitzer]{Eckhard Hitzer}
\address{%
Department of Applied Physics, 
University of Fukui,
Bunkyo 3-9-1, 910-8507 Fukui,
Japan}
\email{hitzer@mech.fukui-u.ac.jp}

\author[R.~Ab\l amowicz]{Rafa\l \ Ab\l amowicz}
\address{%
Department of Mathematics, Box 5054, Tennessee Technological University,
Cookeville, TN 38505, USA}
\email{rablamowicz@tntech.edu}

%\thanks{I thank my family and FTHD organizer S.L. Eriksson.}

\subjclass{Primary 15A66; Secondary 11E88, 42A38, 30G35}

\keywords{Roots of $-1$, Clifford (geometric) algebra, 
Fourier transformation, pseudo scalar}

\date{October 17, 2008}
\dedicatory{Soli Deo Gloria.}
%%% ----------------------------------------------------------------------

\begin{abstract}
It is known that Clifford (geometric) algebra offers a geometric 
interpretation for square roots of $-1$ in the form of blades that 
square to minus $1$. 
This extends to a geometric interpretation of quaternions as the 
side face bivectors of a unit cube. Research has been done~\cite{SJS:Biqroots} 
on the 
biquaternion roots of $-1$, abandoning the restriction 
to blades. Biquaternions are isomorphic to the 
Clifford (geometric) algebra $\G_{3}$ of
$\R^3$. All these roots of 
$-1$ find immediate applications in the construction of new types of 
geometric Clifford Fourier transformations.

We now extend this research to general algebras $\G_{p,q}$. We 
fully derive the geometric roots of $-1$ for the Clifford (geometric) 
algebras with $p+q \leq 4$.
\end{abstract}

% 42A38 Fourier and Fourier-Stieltjes transforms and other transforms of Fourier type
% 11R52 Quaternion and other division algebras: arithmetic, zeta functions
% 15A66 Clifford algebras, spinors 
% 11E88 Quadratic spaces; Clifford algebras
% 30G35 Functions of hypercomplex variables and generalized variables

\maketitle

\section{Introduction}

The British mathematician W.K. Clifford created his \textit{geometric algebras}\footnote{%
In his original publication~\cite{WKC:AppGrass} Clifford first used the term
\textit{geometric algebras}. Subsequently in mathematics the new term 
\textit{Clifford algebras}~\cite{PL:CAaS} has become the proper mathematical term. 
For emphasizing the \textit{geometric} nature of the algebra, some researchers
continue~\cite{DH:NF1,HS:CAtoGC,GS:ncHCFT} to use the original term geometric algebra(s).
} in 1878
inspired by the works of Hamilton on quaternions and by Grassmann's exterior algebra. 
Grassmann invented the antisymmetric outer product of vectors, that regards the oriented parallelogram
area spanned by two vectors as a new type of number, commonly called bivector. The bivector
represents its own plane, because outer products with vectors in the plane vanish. 
In three dimensions the outer product of three linearly independent vectors defines a so-called
trivector with the magnitude of the volume of the parallelepiped spanned by the vectors. 
Its orientation (sign) depends on the handedness of the three vectors. 

In the Clifford algebra~\cite{DH:NF1} of $\R^3$
the three bivector side faces of a unit cube 
$\{\vect{e}_1\vect{e}_2, \vect{e}_2\vect{e}_3, \vect{e}_3\vect{e}_1\}$  
oriented along the three coordinate directions 
$\{\vect{e}_1, \vect{e}_2, \vect{e}_3\}$
correspond to the three quaternion units $\qi$, $\qj$, and $\qk$. Like quaternions,
these three bivectors square to minus one and generate the rotations in their
respective planes. 

Beyond that Clifford algebra allows to extend complex numbers 
to higher dimensions~\cite{HS:CAtoGC,BDS:CA}
and systematically generalize our knowledge of complex numbers, holomorphic functions
and quaternions. It has found rich applications in symbolic computation, physics, robotics, 
computer graphics, etc.~\cite{GS:ncHCFT,TB:thesis,MF:thesis,HL:IAGR}. 
Since bivectors and trivectors in the Clifford algebras of Euclidean
vector spaces square to minus one, we can use them to create new geometric kernels for
Fourier transformations. This leads to a large variety of new Fourier transformations, 
which all deserve to be studied in their own 
right~\cite{LMQ:CAFT94,AM:CAFT96,TQ:PWT,ES:CFTonVF,HM:CFTUP,MHAV:WQFT,EM:CFaUP,
HM:CFToMVF,EH:QFTgen,MHA:2DCliffWinFT,SBS:FastCmQFT,GS:ncHCFT}. 

We will treat both Euclidean (positive definite metric) and non-Euclidean (indefinite metric) vector
spaces. We know from Einstein's special theory of relativity that non-Euclidean
vector spaces are of fundamental importance in nature~\cite{DH:STA}. Therefore this paper is about finding square roots of $-1$ in a non-degenerate Clifford algebra $\G_{p,q}$.

\section{Clifford (geometric) algebras}

The associative geometric product of two vectors $\vect{a},\vect{b} \in \R^{p,q}, p+q=n$
is defined as the sum of their symmetric inner product (scalar) and their antisymmetric outer
product (bivector)
\be\label{eq:gp}
  \vect{a}\vect{b} = \vect{a}\cdot\vect{b}  + \vect{a}\wedge\vect{b}.
\ee
We define~\cite{PL:CAaS} a real Clifford algebra $\G_{p,q}$ as the linear space of all elements generated by
the associative (and distributive) bilinear geometric product of vectors of an inner product
vector space $\R^{p,q}, p+q=n$ over the field of reals $\R$. A Clifford algebra 
includes the field of reals $\R$ and the vector space $\R^{p,q}$ as grade zero and grade one 
elements, respectively. 

Clifford algebras in one, two and three dimensions have the following basis blades of grade 0 (scalars), 
grade 1 (vectors), grade 2 (bivectors) and grade 3 (trivectors)
\be\label{eq:blades}
   \{1, \vect{e}_1, \vect{e}_2, \vect{e}_3, \bv{e}_{23}, \bv{e}_{31}, \bv{e}_{12}, \bv{e}_{123} \},
\ee
where we use abbreviations 
$\bv{e}_{12}=\vect{e}_1\vect{e}_2, \bv{e}_{23}=\vect{e}_2\vect{e}_3, \bv{e}_{31}=\vect{e}_3\vect{e}_1, \bv{e}_{123}=\vect{e}_1\vect{e}_2\vect{e}_3$.
Every multivector can be expanded in terms of these basis blades with real coefficients. We give examples for 
$M \in \G_{p,q}$, $n=p+q=1,2,3$:
\begin{align}
  &M = \alpha + \beta \vect{e}_1, 
  \label{eq:bladeex1}\\
  &M^{\prime} = \alpha + b_1 \vect{e}_1 + b_2 \vect{e}_2 + \beta \bv{e}_{12},
  \label{eq:bladeex2}\\
  &M^{\prime\prime} = \alpha 
  + b_1 \vect{e}_1 + b_2 \vect{e}_2 + b_3 \vect{e}_3
  + c_{1} \bv{e}_{23} + c_{2} \bv{e}_{31} + c_{3} \bv{e}_{12} 
  + \beta \bv{e}_{123} .
  \label{eq:bladeex3}
\end{align}

The general notation for the quadratic form of basis vectors in $\R^{p,q}$ is:
\begin{align}
\vec{e}_k^{\,2} = \varepsilon_k = 
\left\{ 
  \begin{array}{cl}
  +1 & \mbox{ for } 1\leq k \leq p,  \\
  -1 & \mbox{ for } p+1\leq k \leq p+q=n  
  \end{array} 
\right. .
\end{align}
We therefore always have $\vec{e}_k^{\,4} = \varepsilon_k^2 = 1$, and we abbreviate 
$\G_p = \G_{p,0}$. We follow the convention that inner and outer products
have priority over the geometric product, which saves writing a number of brackets. 
Therefore, $\vec{a}\cdot \vec{b} \, \vec{c}$ equals $(\vec{a}\cdot \vec{b}\,) \vec{c}$ and not
$\vec{a} \cdot (\vec{b} \vec{c})$, etc.

We will frequently use the following basic formulas of Clifford algebra in the rest of this work. The symmetric part of the geometric product of any two vectors $\vec{a}$, $\vec{b}$ is the inner product (contraction, scalar product)
\begin{equation}
  \frac{1}{2}(\vec{a}\vec{b} + \vec{b}\vec{a}) 
  = \vec{a}\cdot\vec{b}
  = \langle \vec{a}\vec{b} \rangle_0
  = \langle \vec{a}\vec{b} \rangle.
\end{equation}
Likewise the inner product  (contraction, scalar product) of any two bivectors $\bvc$, $\bvcp$ is symmetric
\begin{equation}
  \frac{1}{2}(\bvc \,\bvcp + \bvcp \bvc) 
  = \bvc \cdot \bvcp
  = \langle \bvc\, \bvcp \rangle_0
  = \langle \bvc\, \bvcp \rangle.
\end{equation}
The antisymmetric part of the geometric product of any two vectors $\vec{a}$, $\vec{b}$ is the outer product (bivector)
\begin{equation}
  \frac{1}{2}(\vec{a}\vec{b} - \vec{b}\vec{a}) 
  = \vec{a}\wedge\vec{b}
  = \langle \vec{a}\vec{b} \rangle_2.
\end{equation}
The inner product (left contraction) of a vector $\vec{a}$ with a bivector $\bvc$ is antisymmetric
\begin{equation}
  \vec{a}\cdot\bvc
  = \frac{1}{2}(\vec{a}\bvc - \bvc \vec{a})
  = \langle \vec{a}\bvc \rangle_1.
\end{equation}
Let $I_n=\Pi_{k=1}^n \vec{e}_k$ be the unit oriented pseudoscalar of $\G_{p,q}, n=p+q$. Let $A_r, B_s \in \G_{p,q}$ be two blades of grade $r$ and $s$, respectively. Then we have the following general rules~\cite{HL:IAGR}. The inner product (left contraction~\cite{LD:InP}) is related to the outer product by\footnote{In order to avoid a discussion of deviating definitions of the inner product for $r=0$ or $s=0$,
we exclude scalars in \eqref{eq:oduali}, but depending on the definition of the
inner product (or contraction), a single general formula for all grades exists. For example for the left contraction~\cite{LD:InP} $A,B \in \G_{p,q}, A\rfloor B = \sum_{r,s}\langle \langle A \rangle_r \langle B \rangle_s\rangle_{s-r} $ we have the two formulas
$(A\wedge B)I_n = A\rfloor (BI_n)$ and $(A\rfloor B)I_n = A\wedge (BI_n)$. }
\begin{equation}
\label{eq:idualo}
  (A_r \cdot B_s)I_n = A_r \wedge (B_s I_n), \mbox{ if } r \leq s;
\end{equation}
\begin{equation}
\label{eq:oduali}
  (A_r \wedge B_s)I_n = A_r \cdot (B_s I_n), \mbox{ if } r + s \leq n, \quad r,s>0.
\end{equation}
Two blades $A_r, B_s \in \G_{p,q}$ are called orthogonal iff their inner product is zero
\begin{equation}
  A_r \perp B_s \quad \Longleftrightarrow \quad A_r \cdot B_s = 0.
\end{equation}
With \eqref{eq:idualo} follows that for $r \leq s$
\begin{equation}
  A_r \perp B_s 
  \quad \Longleftrightarrow \quad A_r \wedge (B_sI_n) = 0
  \quad \Longleftrightarrow \quad A_r \wedge \widetilde{B_s} = 0,
\end{equation}
where  $\widetilde{B_s} = B_s I_n^{-1}$ is the \textit{dual} of $B_s$, with $I_n^{-1} = \pm I_n$. Likewise \eqref{eq:oduali} shows that for $r + s \leq n, \,\,\, r,s>0$
\begin{equation}
\label{eq:perpdual}
  A_r \perp \widetilde{B_s} 
  \quad \Longleftrightarrow \quad A_r \wedge B_s = 0.
\end{equation}
\begin{example}
\label{ex:1}
Let $\vec{b},\bvc \in \G_{p,q}, p+q=3$ be a vector $\vec{b}$ and a bivector $\bvc$  with vanishing outer product. Then by \eqref{eq:perpdual} the dual vector $\vec{c} = \widetilde{\bvc}$ is always perpendicular to $\vec{b}$
independent of the signature of the underlying vector space $\mathbb{R}^{p,q}, p+q=3$,
\begin{equation}
  \vec{b} \wedge \bvc = 0 \quad \Longleftrightarrow \quad 
  \vec{b} \cdot \vec{c} = 0
  \quad \Longleftrightarrow \quad 
  \vec{b} \perp \vec{c}.
\end{equation}
\end{example}

\section{Geometric multivector square roots of $-1$}

\begin{defn}[Geometric root of $-1$]
A geometric multivector square root (geometric root) of $-1$ is a multivector $A\in \G_{p,q}$ with
\begin{align}
  A^2 = AA = -1.
\end{align}
\end{defn}
An immediate application of this definition is the generalization of the famous
Euler formula to geometric roots $A$
\be
  e^{\varphi A} = \cos \varphi + A \sin \varphi .
\ee
For example, Lounesto considers e.g. $\cos \varphi + \bv{e}_{12}\sin \varphi$ in $\G_2$ in~\cite{PL:CAaS} on page 29. 
\begin{thm}
Every multivector square root $A$ of $-1$ is subject to $n+1=p+q+1$ 
grade-wise constraints:
\begin{align}
  \label{eq:re}
  A^2 = \langle AA \rangle = -1,
\end{align}
and
\begin{align}
  \label{eq:cond}
  \langle AA \rangle_k = 0, \qquad 1 \leq  k \leq n,
\end{align}
where $\langle AA\rangle_k $ denotes the $k$-th vector part of $AA$,
and $\langle AA\rangle = \langle AA\rangle_0 $.
\end{thm}

We point out that $\langle AA \rangle$ is identical to the scalar product 
$A \ast A$ of~\cite{HS:CAtoGC}. 
In the following we call the scalar equation \eqref{eq:re} the \textit{root equation} of
$\G_{p,q}$ and \eqref{eq:cond} the \textit{constraints}. Depending on the value of $k$,
each $k$-vector constraint represents $\binom{n}{k}$
scalar equations. We will sometimes conveniently split up a $k$-vector constraint equation
and still call the resulting partial equations \textit{constraints}.

\section{Case $n=1$}

We have two algebras $\G_1$ and $\G_{0,1}$. There is only one basis vector $\vec{e}_1$ with square
$\vec{e}_1^2=\varepsilon_1$. The two Clifford algebras are two dimensional with general elements (multivectors)
\be
   \alpha + \beta \vec{e}_1, \qquad \alpha, \beta \in \R.
\ee
The square of such a multivector is
\be
  (\alpha + \beta \vec{e}_1)^2 = \alpha^2 + \varepsilon_1 \beta^2 + 2\alpha \beta \vec{e}_1 = -1,
\ee
which has the scalar part (root equation)
\be
  \label{eq:scal-1}
  \alpha^2 + \varepsilon_1 \beta^2 = -1,
\ee
and the vector part (constraint)
\be
  \label{eq:vec0}
  2\alpha \beta \vec{e}_1 = 0.
\ee
We see that the left hand side of \eqref{eq:scal-1} is always greater or equal to zero if
$\varepsilon_1=+1$. Therefore $\G_1$  has no multivector square roots of $-1$. 
The vector part \eqref{eq:vec0} is zero if and only if
\be
  \label{eq:vec0a}
  \alpha= 0 \quad \mbox{or} \quad \beta = 0.
\ee
If we try for $\G_{0,1}$ and let $\alpha= 0$, we get  from \eqref{eq:scal-1} the root equation
\be
  \label{eq:scal-1a}
  -1 \beta^2 = -1 \Leftrightarrow \beta^2 = 1 \Leftrightarrow \beta = \pm 1. 
\ee
If we try for $\G_{0,1}$ and let $\beta= 0$, we get  from \eqref{eq:scal-1}
\be
  \label{eq:scal-1b}
  \alpha^2 = -1,
\ee
which is impossible for $\alpha \in \R$. 
Therefore, when  $n=1$, the only geometric roots of $-1$ exist in $\G_{0,1}$
as
\be
  A = \pm \vec{e}_1. 
\ee

\section{Case $n=2$}

We have three central algebras $\G_2$, $\G_{1,1}$ and $\G_{0,2}$. There are two basis vectors $e_k, k\in\{1,2\}$ with square
$e_k^2=\varepsilon_k$. The three Clifford algebras are four dimensional with general elements
\be
   \alpha + \vec{b} +\beta \bv{e}_{12}, \,\,\, 
   \vec{b}=b_1\vec{e}_1+b_2\vec{e}_2, \quad 
   \alpha, b_1, b_2, \beta \in \R, \,\,\,
   \vec{b}\in\R^{p,q}, \,\,\,
   p+q=2.
\ee
The square of such a multivector is
\be
  (\alpha + \vec{b} +\beta \bv{e}_{12})^2 
  = \alpha^2 + \vec{b}^2+ \beta^2 \underbrace{\bv{e}_{12}^2}_{=-\varepsilon_1\varepsilon_2} 
    + 2\alpha\vec{b}  + 2\alpha \beta \bv{e}_{12} + 2\beta \underbrace{(\vec{b}\wedge \bv{e}_{12})}_{=0}= -1,
\ee
which has the scalar part (root equation),
\be
  \label{eq:2scal-1}
  \alpha^2 + \vec{b}^2 -\beta^2 \varepsilon_1\varepsilon_2 = -1,
\ee
two constraints for the vector part 
\be
  \label{eq:2vec0}
  2\alpha\vec{b} = 0,
\ee
and the bivector part
\be
  \label{eq:2biv0}
  2\alpha \beta \bv{e}_{12} = 0.
\ee

\subsection{Case $n=2$, $\alpha=0$}
Equations \eqref{eq:2vec0} and \eqref{eq:2biv0} are now always fulfilled
by any $\vec{b}$ and $\beta$. 
From \eqref{eq:2scal-1} it follows that
\be
  \label{eq:2scal-1b}
  \vec{b}^2 -\beta^2 \varepsilon_1\varepsilon_2 
  =b_1^2 \varepsilon_1 + b_2^2 \varepsilon_2 -\beta^2 \varepsilon_1\varepsilon_2= -1.
\ee
Multiplying each side of (24) by 
$\varepsilon_1\varepsilon_2$ gives the following \textit{root equation}:
\be
  \label{eq:2rteq}
  \beta^2 = b_1^2 \varepsilon_2 + b_2^2 \varepsilon_1 +\varepsilon_1\varepsilon_2
  = \left\{ 
  \begin{array}{cl}
   \phantom{-}b_1^2 + b_2^2 +1  & \mbox{ for } \G_2,\\
  -b_1^2 + b_2^2 -1 & \mbox{ for } \G_{1,1},\\
  -b_1^2 - b_2^2 +1 & \mbox{ for } \G_{0,2}.
  \end{array} 
\right. 
\ee
In $\G_2$ this includes, for $b_1=b_2=0$, the solution $A=\pm \bv{e}_{12}$, which also appears in~\cite{PL:CAaS} on page 29.

\subsection{Case $n=2$, $\alpha \neq 0 $}
If $\alpha \neq 0 $, then, according to \eqref{eq:2vec0} and \eqref{eq:2biv0}, we have  
\be
  \vec{b} = 0 \quad \mbox{and} \quad \beta = 0.
\ee
Inserting this in \eqref{eq:2scal-1} gives
\be
  \alpha^2 = -1, \qquad \alpha \in \R\setminus\{0\},
\ee
which has no solution. Therefore, the root equation \eqref{eq:2rteq} describes already all possible solutions.

\section{Case $n=3$}

We have four algebras $\G_3,\,\G_{2,1},\,\G_{1,2},$ and $\G_{0,3}$ with a non-trivial center  spanned by the identity element $1$ and the unit pseudoscalar $\bv{e}_{123}.$
There are three basis vectors $\vec{e}_k, k\in\{1,2,3\}$, with squares
$\vec{e}_k^2=\varepsilon_k$. The four Clifford algebras are eight dimensional with general elements 
\be
   \alpha + \vec{b} + \underline{c} + \beta \bv{e}_{123}, 
   \quad 
   \alpha, \beta \in \R,
   \quad 
   \vec{b}=b_1\vec{e}_1+b_2\vec{e}_2+b_3\vec{e}_3 \in\R^{p,q},\;p+q=3,  
\ee
with
\be
   \underline{c} = c_1 \bv{e}_{23}+ c_2 \bv{e}_{31}+c_3 \bv{e}_{12}\in \bigwedge^2 \R^{p,q},
   \quad c_1,c_2,c_3 \in \R.   
\ee
Setting the square of such a multivector to $-1$ gives
%\begin{align}
%  (\alpha &+ \vec{b} + \underline{c} + \beta \bv{e}_{123})^2 
%  \nonumber \\
%  = &\,\,\alpha^2 + \vec{b}^2+ \underline{c}^2+ \beta^2 
%    \hspace*{-3mm}\underbrace{\bv{e}_{123}^2}_{=-\varepsilon_1\varepsilon_2\varepsilon_3} 
%    \hspace*{-3mm} + 2\alpha\vec{b}+ 2\alpha \underline{c}  + 2\alpha \beta \bv{e}_{123}  
%  \nonumber \\
%    &+ \underbrace{\vec{b}\underline{c} + \underline{c}\vec{b}}_{=2\vec{b}\wedge\underline{c}}   
%    + 2\beta \vec{b}\bv{e}_{123} + 2\beta \underline{c} \bv{e}_{123} = -1.
%\end{align}
\begin{align}
  (\alpha + \vec{b} + \underline{c} + \beta \bv{e}_{123})^2 &= 
   \alpha^2 + \vec{b}^2+ \underline{c}^2+ \beta^2 
    \hspace*{-3mm}\underbrace{\bv{e}_{123}^2}_{= -\varepsilon_1\varepsilon_2\varepsilon_3} 
    \hspace*{-3mm} + 2\alpha\vec{b}+ 2\alpha \underline{c}  + 2\alpha \beta \bv{e}_{123}  \notag\\
    &\phantom{=}\hspace*{0ex} + \underbrace{\vec{b}\underline{c} + \underline{c}\vec{b}}_{= 2\vec{b}\wedge\underline{c}}   
    + 2\beta \vec{b}\bv{e}_{123} + 2\beta \underline{c} \bv{e}_{123} = -1.
\end{align}
Grade-wise this results in the following set of constraints: 
For the scalar part (root equation)
\be
  \label{eq:3scal-1}
  \alpha^2 + \vec{b}^2+ \underline{c}^2 - \beta^2\varepsilon_1\varepsilon_2\varepsilon_3 = -1,
\ee
for the vector part, 
\be
  \label{eq:3vec0}
  \alpha\vec{b} +  \beta \underline{c} \,\bv{e}_{123} = 0,
\ee
for the bivector part
\be
  \label{eq:3biv0}
  \alpha \underline{c} + \beta \vec{b}\,\bv{e}_{123} = 0,
\ee
and for the trivector part
\be
  \label{eq:3triv0}
  \alpha \beta \bv{e}_{123} + \vec{b}\wedge\underline{c} 
  = (\alpha \beta + b_1c_1+b_2c_2+b_3c_3)\,\bv{e}_{123} 
  = 0.
\ee

\subsection{Case $n=3$, $\alpha=0$}

For $\alpha=0$, the four equations \eqref{eq:3scal-1} to \eqref{eq:3triv0} simplify to
the root equation
\be
  \label{eq:3scal-1a}
  \vec{b}^2+ \underline{c}^2 - \beta^2\varepsilon_1\varepsilon_2\varepsilon_3 = -1,
\ee
and the three constraints
\be
  \label{eq:3notscal0}
  \beta \underline{c} \,\bv{e}_{123} 
  = \beta \vec{b}\,\bv{e}_{123} 
  = \vec{b}\wedge\underline{c} = (b_1c_1+b_2c_2+b_3c_3)\,\bv{e}_{123} 
  = 0.
\ee
The expression $\vec{b}\wedge\underline{c}=0$ means that $\vec{b}$ is in the plane 
defined by the bivector $\underline{c}$, which can also be written as
\be
  \label{eq:3notscal0a}
  \vec{b}\underline{c} = \vec{b} \cdot \underline{c} \,.
\ee
In three dimensions the bivector $\underline{c}$ can also be represented
by its dual vector (perpendicular to the plane defined by $\underline{c}$)
\begin{align}
  \vec{c}
  = \underline{c}\bv{e}_{123}^{-1}
  = \varepsilon_1 c_1 \vec{e}_1+\varepsilon_2 c_2\vec{e}_2+\varepsilon_3 c_3\vec{e}_3
  = \left\{
  \begin{array}{ll}
  c_1 \vec{e}_1+c_2\vec{e}_2+c_3\vec{e}_3 & \text{for } \G_3, \\
  c_1 \vec{e}_1+c_2\vec{e}_2-c_3\vec{e}_3 & \text{for } \G_{2,1}, \\
  c_1 \vec{e}_1-c_2\vec{e}_2-c_3\vec{e}_3 & \text{for } \G_{1,2}, \\
  -c_1 \vec{e}_1-c_2\vec{e}_2-c_3\vec{e}_3\hspace*{-3mm} & \text{for } \G_{0,3}, 
  \end{array}
  \right.
\end{align}
where we used $\varepsilon_k = \pm 1, \, k \in \{ 1,2,3 \}$ and hence  
$\varepsilon_k^{-1} = \varepsilon_k$.
Therefore, \textit{independent} of the signature of the quadratic form, we have the following constraint
\be 
  \label{eq:3triv0a}
  (\vec{b} \wedge  \bvc) \, \bv{e}_{123}^{-1}= 
  \vec{b} \cdot  \vec{c} = b_1c_1+b_2c_2+b_3c_3 = 0 
  \,\,\, \Longleftrightarrow \,\,\, 
  \vec{b} \vec{c} = \vec{b} \wedge \vec{c} \,,
\ee
i.e., $\vec{b} \perp \vec{c}$, which should be compared with Example \ref{ex:1}.

\subsubsection{Case $n=3$, $\alpha=0$, $\beta=0$\label{sc:3al0be0}}
The constraints are now given by $\alpha=0$, $\beta=0$, 
and \eqref{eq:3notscal0a} or \eqref{eq:3triv0a}.
For $\alpha = \beta=0 $ the root equation \eqref{eq:3scal-1a} further simplifies to
\begin{align}
  \label{eq:3scal-1b}
  -1 &= \vec{b}^2+ \underline{c}^2 
  = b_1^2\varepsilon_1 +  b_2^2\varepsilon_2 + b_3^2\varepsilon_3
    - c_1^2 \varepsilon_2\varepsilon_3 - c_2^2 \varepsilon_3\varepsilon_1 - c_3^2 \varepsilon_1\varepsilon_2 
  \nonumber \\
  &=
  \left\{
    \begin{array}{ll}
      \vec{b}^2-\vec{c}^{\,\,2}  
      & 
      = \left\{\begin{array}{ll}
      b_1^2 +  b_2^2 + b_3^2 - (c_1^2 + c_2^2 + c_3^2) & \text{ for } \G_{3},\\
      b_1^2 -  b_2^2 - b_3^2 - (c_1^2 - c_2^2 - c_3^2) & \text{ for } \G_{1,2},
      \end{array} \right. \vspace*{2mm}
      \\
      \vec{b}^2+\vec{c}^{\,\,2} 
      & 
      = \left\{\begin{array}{ll}
      b_1^2 +  b_2^2 - b_3^2 + (c_1^2 + c_2^2 - c_3^2)     & \text{ for } \G_{2,1},\\
      - (b_1^2 +  b_2^2 + b_3^2) - (c_1^2 + c_2^2 + c_3^2) & \text{ for } \G_{0,3},
      \end{array} \right.
    \end{array}
  \right. 
\end{align}
and \eqref{eq:3triv0a}.
We now explain the geometric interpretation of the root equations \eqref{eq:3scal-1b}.

For $\G_{3}$ equation \eqref{eq:3scal-1b} means the perpendicular vectors $\vec{b}$ and 
$\vec{c}$ define a quadric (a 6D hyperboloid) in $\R^6$.
Or in other words, for a given vector $\vec{b}$ the bivectors $\underline{c}$ that
lead to geometric roots of $-1$ are defined by all radial vectors $\vec{c}$ 
of a circle in a plane perpendicular to $\vec{b}$ 
with radius $|\vec{c}| = \sqrt{1+\vec{b}^2}$.

For $\G_{2,1}$ and $\G_{1,2}$ the respective equations in \eqref{eq:3scal-1b} define quadrics of possible solutions in the space $\R^6 = \R^3 \oplus \R^3$ of vectors $\vec{b}$
and dual vectors $\vec{c}$. 

For $\G_{0,3}$ the respective equation \eqref{eq:3scal-1b} defines again a quadric of possible solutions in the space $\R^6 = \R^3 \perp \R^3$ of vectors $\vec{b}$
and dual vectors $\vec{c}$ 
perpendicular to $\vec{b}$. 
The quadric in question can be pictured as a unit sphere in $\R^6$. Geometric roots
of $-1$ exist only for vectors $\vec{b}$ with $|\vec{b}|\leq 1$. The possible dual
vectors $\vec{c}$ are the radial vectors of a circle defined by the intersection of
a plane (with distance $\vec{b}$ from the origin) with the $\mathbb{R}^3$ unit sphere (centered
at the origin).

\subsubsection{Case $n=3$, $\alpha=0$, $\beta \neq 0$\label{sc:3a0bn0}}
For $\beta \neq 0$, equation \eqref{eq:3notscal0} simplifies to
\be
  \label{eq:3notscal0b}
  \underline{c}=0, \qquad
  \vec{b} = 0,
\ee
while the root equation \eqref{eq:3scal-1a} gives
\be
  \label{eq:3scal-1g}
  \beta^2\varepsilon_1\varepsilon_2\varepsilon_3 = 1.
\ee
Because for real $\beta\in\R\setminus\{0\}$ the square $\beta^2>0$ is 
always positive, equation \eqref{eq:3scal-1g} includes two constraints
\be
  \beta = \pm 1 \quad \mbox{and} \quad \varepsilon_1\varepsilon_2\varepsilon_3 = 1.
\ee
This is only possible in $\G_{3}$ and $\G_{1,2}$, 
but not in $\G_{2,1}$ and $\G_{0,3}$. So for $\G_{3}$ and $\G_{1,2}$ we get
the geometric trivector roots of $-1$ as
\be
  A = \pm \,\bv{e}_{123}.
\ee

\subsection{Case $n=3$, $\alpha\neq 0$}
We will see that no more geometric roots of $-1$ arise for the case $\alpha\neq 0$.
To prove this is not trivial as we will see in the following. 

For $\alpha\neq 0$ we get from \eqref{eq:3vec0}
\be
  \label{eq:3vec0a}
  \vec{b} = - \frac{\beta}{\alpha} \,\,\underline{c} \,\bv{e}_{123},
\ee
from \eqref{eq:3biv0} 
\be
  \label{eq:3biv0a}
  \underline{c} = - \frac{\beta}{\alpha} \,\,\vec{b}\,\bv{e}_{123},
\ee
and from \eqref{eq:3triv0} 
\be
  \label{eq:3triv0b}
  \vec{b}\wedge\underline{c} 
  =(b_1c_1+b_2c_2+b_3c_3)\,\bv{e}_{123} 
  = -\alpha \beta \bv{e}_{123}. 
\ee
By squaring both sides of equations \eqref{eq:3vec0a} and \eqref{eq:3biv0a} we obtain
\be
  \label{eq:3vec0a1}
  \vec{b}^{\,2} = \frac{\beta^2}{\alpha^2} \,\,\underline{c}^2 \,\bv{e}_{123}^2,
\ee
and
\be
  \label{eq:3biv0a1}
  \underline{c}^2 = \frac{\beta^2}{\alpha^2} \,\,\vec{b}^{\,2}\,\bv{e}_{123}^2.
\ee
Inserting \eqref{eq:3biv0a1} into \eqref{eq:3vec0a1} yields
\be
  \label{eq:3vec0a2}
  \vec{b}^{\,2} = \frac{\beta^4}{\alpha^4} \,\,\vec{b}^{\,2}.
\ee
If $\vec{b}^{\,2} \neq 0$, we get from \eqref{eq:3vec0a2} 
that 
\be
  \label{eq:albe}
  \beta^2 = \alpha^2 \quad \mbox{and, therefore,} \quad \beta = \pm \alpha .
\ee

\subsubsection{Case $n=3$, $\alpha\neq 0$, $\beta=0$}
For $\alpha\neq 0$ and $\beta=0$ equations \eqref{eq:3vec0a} and \eqref{eq:3biv0a} 
further simplify to
\be
  \label{eq:3allanot0b0}
  \vec{b} = 0, 
  \quad
  \underline{c} = 0.
\ee
Equation \eqref{eq:3triv0b} is then trivially fulfilled. 
The root equation \eqref{eq:3scal-1} reduces to
\be
  \alpha^2 = -1,
\ee
which cannot be fulfilled for $\alpha\in\R\setminus\{0\}$. Therefore no geometric roots
of $-1$ exist for $\alpha\neq 0$ and $\beta=0$.

\subsubsection{Case $n=3$, $\alpha\neq 0$, $\beta\neq0$}
We now insert \eqref{eq:3biv0a} into \eqref{eq:3triv0b} to get
\begin{align}
  \label{eq:3triv0b1}
  &\vec{b}\wedge \left(-\frac{\beta}{\alpha} \vec{b}\,\bv{e}_{123}\right)  
  = -\alpha \beta \bv{e}_{123}
  \,\,\, \stackrel{\beta \neq 0}{\Longleftrightarrow} \,\,\,
  (\vec{b} \cdot \vec{b}\,) \, \bv{e}_{123} = \alpha^2 \bv{e}_{123}
  \nonumber \\
  &\,\,\, {\Longleftrightarrow} \,\,\,
  \vec{b}^{\,2} = \alpha^2
  \,\,\, \stackrel{\eqref{eq:3biv0a1}}{\Longrightarrow} \,\,\,
  \underline{c}^2 = \beta^2 \bv{e}_{123}^2.
\end{align}
If $\alpha\neq 0$, then also $\alpha^2\neq 0$ and therefore according to \eqref{eq:3triv0b1}
$\vec{b}^{\,2} \neq 0$. According to \eqref{eq:albe} we now have $\beta^2 = \alpha^2 $.
Inserting $\vec{b}^{\,2}$, $\underline{c}^2$ and $\beta^2$ into the root equation \eqref{eq:3scal-1}
gives
\be
  \label{eq:3scal-1c}
  \alpha^2 + \alpha^2+ \alpha^2 \bv{e}_{123}^2 + \alpha^2 \bv{e}_{123}^2
  = 2 \alpha^2 (1+ \bv{e}_{123}^2) = -1 .
\ee
Because $\bv{e}_{123}^2=-1$ for $\G_{3}$ and $\G_{1,2}$, 
and $\bv{e}_{123}^2=+1$ for $\G_{2,1}$ and $\G_{0,3}$, 
we get from \eqref{eq:3scal-1c} the root equations
\be
  \label{eq:0e-1}
  0 = -1 \text{ for } \G_{3} \text{ and } \G_{1,2},
\ee
and 
\be
  \label{eq:4ae-1}
  4 \alpha^2 = -1 \text{ for } \G_{2,1} \text{ and } \G_{0,3}.
\ee
For real $\alpha \neq 0$ both \eqref{eq:0e-1} and \eqref{eq:4ae-1} have no solution.

Therefore the only geometric roots of $-1$ for $n=3$
are the ones found in section \ref{sc:3al0be0} for $\alpha=\beta=0$, 
and in section \ref{sc:3a0bn0} for $\alpha=0$, $\beta\neq0$. 
No geometric roots of $-1$ for $n=3$ exist for $\alpha\neq0$.

The geometric roots of $-1$ of $\G_{p,q}$, $n=p+q\leq 3$ are summarized
in Table~\ref{tb:roots123} on page \pageref{pg:roots123}. We point out, 
that the root equation for $n=2, \alpha = 0$ results from simply inserting the case condition $\alpha=0$ of column two into the general $n=2$ root equation \eqref{eq:2scal-1}. Likewise, the root equation for $n=3, \alpha = \beta = 0$
results from simply inserting the case condition $\alpha=\beta=0$ of column two into the general $n=3$ root equation \eqref{eq:3scal-1}.

\section{Case $n=4$}

We have five central algebras $\G_4$, $\G_{3,1}$, $\G_{2,2}$, $\G_{1,3}$, and $\G_{0,4}$. 
There are four basis vectors $\vec{e}_k, k\in\{1,2,3,4\}$ with square
$\vec{e}_k^{\,2}=\varepsilon_k$, 
$\vec{e}_k^{\,4}=\varepsilon_k^2 = 1$,
$\bv{e}_{123}^2 = \bv{e}_{123}^{-2}$,
and
$\bv{e}_{123}^4 = 1$.
The five Clifford algebras are 16 dimensional with general elements
\begin{align}
   \label{eq:4multiv}
   &\alpha + \vb + \bvc + \beta \, \bv{e}_{123} 
   + (\alphap + \vbp + \bvcp + \betap \, \bv{e}_{123})\vect{e}_4, 
   \nonumber \\ 
   &\alpha, \beta, \alphap, \betap \in \R,
   \nonumber \\ 
   &\vb=b_1\vec{e}_1+b_2\vec{e}_2+b_3\vec{e}_3, \,\,\,
    \vbp=b_1^{\prime}\vec{e}_1+b_2^{\prime}\vec{e}_2+b_3^{\prime}\vec{e}_3 \in\R^{p,q},p+q=3, 
   \nonumber \\ 
   &\bvc = c_1 \bv{e}_{23}+ c_2 \bv{e}_{31}+c_3\bv{e}_{12}, \,\,\,
    \bvcp = c_1^{\prime} \bv{e}_{23}+ c_2^{\prime} \bv{e}_{31}+c_3^{\prime}\bv{e}_{12} 
    \in \bigwedge^2 \R^{p,q}.
\end{align}
Setting the square of such a multivector to $-1$ gives:
\begin{align}
  &[\alpha + \vb + \bvc + \beta \, \bv{e}_{123} 
   + (\alphap + \vbp + \bvcp + \betap \, \bv{e}_{123})\vect{e}_4]^2 
  \nonumber \\
  &= (\alpha + \vb + \bvc + \beta \, \bv{e}_{123})^2
  + (\alphap + \vbp + \bvcp + \betap \, \bv{e}_{123})\vect{e}_4(\alphap + \vbp + \bvcp + \betap \, \bv{e}_{123})\vect{e}_4
  \nonumber \\
  &+ (\alpha + \vb + \bvc + \beta \, \bv{e}_{123})(\alphap + \vbp + \bvcp + \betap \, \bv{e}_{123})\vect{e}_4
  \nonumber \\  &
  + (\alphap + \vbp + \bvcp + \betap \, \bv{e}_{123})\vect{e}_4(\alpha + \vb + \bvc + \beta \, \bv{e}_{123})
  \nonumber \\
  &= (\alpha + \vb + \bvc + \beta \, \bv{e}_{123})^2
  + (\alphap + \vbp + \bvcp + \betap \, \bv{e}_{123})(\alphap - \vbp + \bvcp - \betap \, \bv{e}_{123})\varepsilon_4
  \nonumber \\
  &+ (\alpha + \vb + \bvc + \beta \, \bv{e}_{123})(\alphap + \vbp + \bvcp + \betap \, \bv{e}_{123})\vect{e}_4
  \nonumber \\  &
  + (\alphap + \vbp + \bvcp + \betap \, \bv{e}_{123})(\alpha - \vb + \bvc - \beta \, \bv{e}_{123})\vect{e}_4
  \nonumber \\
  &= -1.
\end{align}
We therefore get
%\begin{align}
\begin{multline}
  \label{eq:n4-1}
  (\alpha + \vb + \bvc + \beta \, \bv{e}_{123})^2 \\
 + (\alphap + \vbp + \bvcp + \betap \, \bv{e}_{123})(\alphap - \vbp + \bvcp - \betap \, \bv{e}_{123})\varepsilon_4
  = -1,
%\end{align}
\end{multline}
and
%\begin{align}
\begin{multline}
  \label{eq:4e4part}
(\alpha + \vb + \bvc + \beta \, \bv{e}_{123})(\alphap + \vbp + \bvcp + \betap \, \bv{e}_{123})\vect{e}_4 \\
  +(\alphap + \vbp + \bvcp + \betap \, \bv{e}_{123})(\alpha - \vb + \bvc - \beta \, \bv{e}_{123})\vect{e}_4 = 0.
%\end{align}
\end{multline}
Multiplying out \eqref{eq:n4-1} gives
\begin{align}
  \label{eq:n4-1a}
  &\alpha^2 + \vb^2 + \bvc^2 + \beta^2 \, \bv{e}_{123}^2
  + \varepsilon_4 \alphap^2 
  - \varepsilon_4 \vbpsq + \varepsilon_4\bvcp^2 - \varepsilon_4\betap^2 \, \bv{e}_{123}^2
  \nonumber \\
  &+ 2 \alpha \vb + 2 \alpha \bvc + 2 \alpha \beta \, \bv{e}_{123}
  + \underbrace{\vb \bvc + \bvc \vb}_{2\vb \wedge \bvc} 
  + \underbrace{\beta \vb \, \bv{e}_{123} + \beta \, \bv{e}_{123} \vb }_{2\beta \vb \, \bv{e}_{123}}
  + \underbrace{\beta \bvc \, \bv{e}_{123} + \beta \, \bv{e}_{123} \bvc}_{2 \beta \bvc \, \bv{e}_{123}}
  \nonumber \\
  &+ [ \underbrace{- \alphap \vbp + \alphap \vbp }_{0} 
    + 2\alphap \bvcp 
    \underbrace{-\alphap \betap \, \bv{e}_{123} +\alphap \betap \, \bv{e}_{123}}_{0}
    \underbrace{-\betap \vbp \, \bv{e}_{123} -\betap \, \bv{e}_{123} \vbp}_{-2\betap \vbp \, \bv{e}_{123}}
  \nonumber \\
  & + \underbrace{\vbp \bvcp - \bvcp \vbp}_{2\vbp \cdot \bvcp}
    \underbrace{-\betap \bvcp \, \bv{e}_{123} + \betap \, \bv{e}_{123} \bvcp}_{0}
  ] \varepsilon_4
  =-1.
\end{align}
This results grade-wise in the following set of equations. 
For the scalar part (\textit{root equation})
\be
  \label{eq:4root}
  \alpha^2 + \vec{b}^2+ \underline{c}^2 + \beta^2\, \bv{e}_{123}^2 
  + \varepsilon_4 \alphap^2 
  - \varepsilon_4 \vbpsq + \varepsilon_4\bvcp^2 - \varepsilon_4\betap^2 \, \bv{e}_{123}^2
  = -1,
\ee
the vector part of the l.h.s. in \eqref{eq:n4-1}
\be
  \label{eq:4_v0}
  \alpha\vec{b} +  \beta \underline{c} \, \bv{e}_{123} + \varepsilon_4\vbp\cdot \bvcp= 0,
\ee
the bivector part of the l.h.s. in \eqref{eq:n4-1}
\be
  \label{eq:4_bv0}
  \alpha \underline{c}  + \varepsilon_4\alphap \bvcp + (\beta \vec{b} - \varepsilon_4\betap \vbp) \, \bv{e}_{123}= 0,
\ee
and the trivector part of the l.h.s. in \eqref{eq:n4-1}
\be
  \label{eq:4_t0}
  \alpha \beta \, \bv{e}_{123} + \vec{b}\wedge\underline{c} 
  = (\alpha \beta + b_1c_1+b_2c_2+b_3c_3)\, \bv{e}_{123} 
  = 0,
\ee
%The $e_4$ equation \eqref{eq:4e4part} results in
After multiplying both sides of equation \eqref{eq:4e4part} by $(\e{4})^{-1}$ we get
\begin{align}
  \label{eq:4e4part_a}
  &\alpha \alphap + \vb \vbp + \bvc \,\bvcp + \beta \betap \, \bv{e}_{123}^2 
  + \alpha \vbp + \alphap \vb + \alpha \bvcp + \alphap \bvc + \alpha \betap \, \bv{e}_{123} + \alphap \beta \, \bv{e}_{123}
  \nonumber\\
  &+ \vb \bvcp + \bvc \vbp + (\beta \vbp + \betap \vb)\, \bv{e}_{123} + \betap \bvc \, \bv{e}_{123} + \beta \bvcp \, \bv{e}_{123}
  \nonumber \\  
  &+ \alpha \alphap - \vbp \vb + \bvcp \bvc - \beta \betap \, \bv{e}_{123}^2 - \alphap \vb + \alpha \vbp 
  + \alphap \bvc + \alpha \bvcp + \alpha \betap \, \bv{e}_{123} - \alphap \beta \, \bv{e}_{123} 
  \nonumber \\
  &+ \vbp \bvc - \bvcp \vb -(\beta \vbp + \betap \vb)\, \bv{e}_{123} + (\betap \bvc - \beta \bvcp)\, \bv{e}_{123} 
  = 0.   
\end{align}
Simplification of \eqref{eq:4e4part_a}, similar to \eqref{eq:n4-1a}, results in
%\begin{align}
%  \label{eq:4e4part_b}
%  &2\alpha \alphap + 2\vb \wedge \vbp +  2\bvc \cdot \bvcp 
%  + 2\alpha \vbp + 2\alpha \bvcp + 2\alphap \bvc + 2\alpha \betap \, \bv{e}_{123} 
%  \nonumber \\  
%  &+ 2 \vbp \wedge \bvc +2\vb \cdot \bvcp + 2 \betap \bvc \, \bv{e}_{123} = 0.
%\end{align}
\begin{multline}
  \label{eq:4e4part_b}
  2\alpha \alphap + 2\vb \wedge \vbp +  2\bvc \cdot \bvcp 
  + 2\alpha \vbp + 2\alpha \bvcp + 2\alphap \bvc + 2\alpha \betap \, \bv{e}_{123} 
  \\  
  + 2 \vbp \wedge \bvc +2\vb \cdot \bvcp + 2 \betap \bvc \, \bv{e}_{123} = 0.
\end{multline}
Grade-wise we get from \eqref{eq:4e4part_b} the scalar part
\begin{align}
  \label{eq:4e4part_s}
  &\alpha \alphap + \bvc \cdot \bvcp = 0,
\end{align}
the vector part
\begin{align}
  \label{eq:4e4part_v}
  &\alpha \vbp + \vb \cdot \bvcp + \betap \bvc \, \bv{e}_{123} = 0,
\end{align}
the bivector part
\begin{align}
  \label{eq:4e4part_bv}
  &\vb \wedge \vbp + \alpha \bvcp + \alphap \bvc = 0,
\end{align}
and the trivector part
\begin{align}
  \label{eq:4e4part_t}
  & \alpha \betap \, \bv{e}_{123} + \vbp \wedge \bvc = 0.
\end{align}
Apart from the actual root equation \eqref{eq:4root} we have therefore the following set of seven constraint equations
\begin{align}
  \label{eq:4e4part_s1}
  &\bvc \cdot \bvcp = - \alpha \alphap,
\end{align}
\be
  \label{eq:4_v1}
  \alpha\vec{b} = - \varepsilon_4\vbp\cdot \bvcp - \beta \underline{c} \, \bv{e}_{123} ,
\ee
\begin{align}
  \label{eq:4e4part_v1}
  \alpha \vbp = - \vb \cdot \bvcp - \betap \bvc \, \bv{e}_{123},
\end{align}
\begin{align}
  \label{eq:4e4part_bv1}
  \alpha \bvcp + \alphap \bvc = \vbp \wedge \vb,
\end{align}
\be
  \label{eq:4_bv1}
  \alpha \underline{c}  + \varepsilon_4\alphap \bvcp = (\varepsilon_4\betap \vbp - \beta \vec{b}) \, \bv{e}_{123},
\ee
\be
  \label{eq:4_t1}
   - \vec{b}\wedge\underline{c} = - \underline{c} \wedge \vec{b}= \alpha \beta \, \bv{e}_{123},
\ee
\begin{align}
  \label{eq:4e4part_t1}
  - \vbp\wedge \bvc = - \bvc \wedge \vbp = \alpha \betap \, \bv{e}_{123}.
\end{align}
The outer products of \eqref{eq:4e4part_bv1} with $\vb$ and $\vbp$ give the following 
useful identities
\begin{align}
  \label{eq:4e4part_bv2}
  \alpha \vb \wedge \bvcp + \alphap \vb \wedge \bvc = 0 
  \stackrel{\eqref{eq:4_t1}}{\Longrightarrow} 
  \alpha \,\vb \wedge \bvcp = \alpha \alphap \beta \, \bv{e}_{123},
\end{align}
\begin{align}
  \label{eq:4e4part_bv3}
  \alpha \vbp \wedge \bvcp + \alphap \vbp \wedge \bvc = 0
   \stackrel{\eqref{eq:4e4part_t1}}{\Longrightarrow} 
  \alpha \,\vbp \wedge \bvcp = \alpha \alphap \betap \, \bv{e}_{123}.
\end{align}
The inner products (left contractions) of \eqref{eq:4_v1} with $\vbp$ and of \eqref{eq:4e4part_v1} with $\vb$
lead to
\be
  \label{eq:4_v2}
  \alpha\vb\cdot \vbp 
  = - \varepsilon_4\underbrace{\vbp\cdot(\vbp\cdot \bvcp)}_{0} 
    - \beta \underbrace{\vbp\cdot(\bvc \, \bv{e}_{123})}_{(\vbp\wedge\bvc) \, \bv{e}_{123}}
  = \beta(- \vbp\wedge\bvc) \, \bv{e}_{123}
  \stackrel{\eqref{eq:4e4part_t1}}{=} \alpha \beta\betap \, \bv{e}_{123}^2,
\ee
\begin{align}
  \label{eq:4e4part_v2}
  \alpha \vb \cdot \vbp 
  = - \underbrace{\vb\cdot(\vb \cdot \bvcp)}_{0} 
    - \betap \underbrace{\vb\cdot(\bvc \, \bv{e}_{123})}_{(\vb\wedge\bvc) \bv{e}_{123}}
  = \betap (-\vb\wedge\bvc) \, \bv{e}_{123}
  \stackrel{\eqref{eq:4_t1}}{=} \alpha \beta\betap \, \bv{e}_{123}^2.
\end{align}

\noindent
We further contract each side of \eqref{eq:4_bv1} from the left with $\bvc$ to obtain
\begin{align}
  \label{eq:4_bv2}
  &\alpha \bvc^2  + \varepsilon_4\alphap \underbrace{\bvc\cdot\bvcp}_{-\alpha \alphap}
  \nonumber \\
  &= \bvc\cdot[(\varepsilon_4\betap \vbp - \beta \vec{b}) \, \bv{e}_{123}]
  \nonumber \\
  &= \varepsilon_4\betap (\bvc\wedge \vbp) \, \bv{e}_{123} - \beta (\bvc\wedge\vec{b})\, \bv{e}_{123}
  \nonumber \\
  &\stackrel{\eqref{eq:4_t1},\eqref{eq:4e4part_t1}}{=} - \varepsilon_4\alpha \betap^2 \, \bv{e}_{123}^2 + \alpha \beta^2 \, \bv{e}_{123}^2,
\end{align}
or, equivalently,
\begin{align}
  \label{eq:4_bv3}
  \alpha \bvc^2  
  = \alpha [\varepsilon_4 \alphap^2 
    - \varepsilon_4 \betap^2 \, \bv{e}_{123}^2 
    + \beta^2 \, \bv{e}_{123}^2].
\end{align}
For $\alpha \neq 0,$ we similarly contract each side of \eqref{eq:4_bv1} from the left with 
$\bvcp$ to obtain
\begin{align}
  \label{eq:4_bv4}
  &\alpha \underbrace{\bvc\cdot\bvcp}_{-\alpha \alphap}  + \varepsilon_4\alphap \bvcp^2
  \nonumber \\
  &= \bvcp\cdot[(\varepsilon_4\betap \vbp - \beta \vec{b}) \, \bv{e}_{123}]
  \nonumber \\
  &= \varepsilon_4\betap (\bvcp\wedge \vbp) \, \bv{e}_{123} - \beta (\bvcp\wedge\vec{b})\, \bv{e}_{123}
  \nonumber \\
  & \stackrel{\eqref{eq:4e4part_bv2},\eqref{eq:4e4part_bv3}}{=} \varepsilon_4\alphap \betap^2 \, \bv{e}_{123}^2 - \alphap \beta^2 \, \bv{e}_{123}^2,
\end{align}
or equivalently ($\varepsilon_4^2 =1$)
\begin{align}
  \label{eq:4_bv5}
  \alphap \bvcp^2
  = \alphap [\varepsilon_4\alpha^2  
    +  \betap^2 \, \bv{e}_{123}^2 -  \varepsilon_4\beta^2 \, \bv{e}_{123}^2].
\end{align}
\noindent
The inner product of \eqref{eq:4_v1} with $\vb$ 
leads to
\begin{align}
  \label{eq:4_v3}
  \alpha\vb^2 
  &= - \varepsilon_4\underbrace{\vb\cdot(\vbp\cdot \bvcp)}_{(\vb\wedge\vbp)\cdot \bvcp} 
    - \beta \underbrace{\vb\cdot(\bvc \, \bv{e}_{123})}_{(\vb\wedge\bvc) \, \bv{e}_{123}}
  \nonumber \\
  &= \varepsilon_4 (\alpha \bvcp + \alphap \bvc )\cdot \bvcp
    + \alpha \beta^2 \, \bv{e}_{123}^2
  \nonumber \\
  &= \varepsilon_4 \alpha \bvcp^2 
     + \varepsilon_4\alphap \underbrace{\bvc \cdot \bvcp}_{-\alpha \alphap}
     +\alpha \beta^2 \, \bv{e}_{123}^2
  \nonumber \\
  &= \alpha [\varepsilon_4 \bvcp^2
     -\varepsilon_4 \alphap^2 + \beta^2 \, \bv{e}_{123}^2]
\end{align}
where we inserted \eqref{eq:4e4part_bv1} and \eqref{eq:4_t1} for the second equality.
Assuming $\alphap \neq 0$, equation \eqref{eq:4_v3} leads with \eqref{eq:4_bv5} to
\begin{align}
  \label{eq:4_v4}
  \alpha\vb^2 
  &= \varepsilon_4 \alpha [\varepsilon_4\alpha^2  
    +  \betap^2 \, \bv{e}_{123}^2 -  \varepsilon_4\beta^2 \, \bv{e}_{123}^2] 
     -\varepsilon_4\alpha \alphap^2 +\alpha \beta^2 \, \bv{e}_{123}^2
  \nonumber \\
  &= \alpha [\alpha^2 -\varepsilon_4 \alphap^2 +  \varepsilon_4 \betap^2 \, \bv{e}_{123}^2] ,
\end{align}
\noindent
The inner product of \eqref{eq:4e4part_v1} with $\vbp$ leads to
\begin{align}
  \label{eq:4e4part_v3}
  &\alpha \vec{b}^{\prime 2}
  = - \underbrace{\vbp\cdot(\vb \cdot \bvcp)}_{(\vbp\wedge\vb)\cdot \bvcp} 
    - \betap \underbrace{\vbp\cdot(\bvc \, \bv{e}_{123})}_{(\vbp\wedge\bvc) \, \bv{e}_{123}}
  \nonumber \\
  &= -(\alpha \bvcp + \alphap \bvc )\cdot \bvcp
    + \alpha \betap^2 \, \bv{e}_{123}^2
  \nonumber \\
  &= -\alpha \bvcp^2 - \alphap \underbrace{\bvc \cdot \bvcp}_{-\alpha \alphap} + \alpha \betap^2 \, \bv{e}_{123}^2
  \nonumber \\
  &= \alpha [-\bvcp^2 + \alphap^2 + \betap^2 \, \bv{e}_{123}^2],
\end{align}
where we inserted \eqref{eq:4e4part_bv1} and \eqref{eq:4e4part_t1} for the second equality.
Assuming $\alphap \neq 0$, equation \eqref{eq:4e4part_v3} leads with \eqref{eq:4_bv5} to
\begin{align}
  \label{eq:4e4part_v4}
  \alpha \vec{b}^{\prime 2} &= -\alpha [\varepsilon_4\alpha^2  
    +  \betap^2 \, \bv{e}_{123}^2 -  \varepsilon_4\beta^2 \, \bv{e}_{123}^2] +\alpha \alphap^2 + \alpha \betap^2 \, \bv{e}_{123}^2
  \notag \\
  &= \alpha [\alphap^2 -\varepsilon_4\alpha^2 +  \varepsilon_4\beta^2 \, \bv{e}_{123}^2] . 
\end{align}
\noindent
Inserting \eqref{eq:4_bv3}, \eqref{eq:4_v3}, and \eqref{eq:4e4part_v3} 
into the root equation \eqref{eq:4root} for $\alpha\neq 0$ we obtain (for all $\alphap$)
\begin{align}
  \label{eq:4root_a}
  &\alpha^2 + \vec{b}^2+ \underline{c}^2 + \beta^2\, \bv{e}_{123}^2 
  + \varepsilon_4 \alphap^2 
  - \varepsilon_4 \vbpsq + \varepsilon_4\bvcp^2 - \varepsilon_4\betap^2 \, \bv{e}_{123}^2
  \nonumber \\
  &= \alpha^2 
     + \varepsilon_4 \bvcp^2 -\varepsilon_4 \alphap^2 + \beta^2 \, \bv{e}_{123}^2
     + \varepsilon_4 \alphap^2 - \varepsilon_4 \betap^2 \, \bv{e}_{123}^2 + \beta^2 \, \bv{e}_{123}^2
     + \beta^2 \, \bv{e}_{123}^2
  \nonumber \\
     &+ \varepsilon_4 \alphap^2
     + \varepsilon_4 \bvcp^2 - \varepsilon_4 \alphap^2 - \varepsilon_4 \betap^2 \, \bv{e}_{123}^2
     + \varepsilon_4 \bvcp^2
     - \varepsilon_4 \betap^2 \, \bv{e}_{123}^2
  \nonumber \\
  &= \alpha^2 +3\alpha^2 -3\alpha^2 +3\varepsilon_4 \bvcp^2 + 3\beta^2 \, \bv{e}_{123}^2 -3 \varepsilon_4\betap^2\, \bv{e}_{123}^2
  \nonumber \\
  &= 4\alpha^2  +3\varepsilon_4 [\bvcp^2 - \varepsilon_4\alpha^2 -\betap^2\, \bv{e}_{123}^2 + \varepsilon_4\beta^2 \, \bv{e}_{123}^2] 
  = -1,
\end{align}
If in addition $\alphap\neq 0$ then with \eqref{eq:4_bv5} we get for the root equation
\begin{align}
  \label{eq:4root_b}
  &\alpha^2 + \vec{b}^2+ \underline{c}^2 + \beta^2\, \bv{e}_{123}^2 
  + \varepsilon_4 \alphap^2 
  - \varepsilon_4 \vbpsq + \varepsilon_4\bvcp^2 - \varepsilon_4\betap^2 \, \bv{e}_{123}^2
  \nonumber \\
  & = 4\alpha^2 + 0 = -1,
\end{align}
Therefore, we have no solution for $\alpha\neq 0$ and $\alphap\neq 0$.

\subsection{$n=4$, $\alpha \neq  0$, $\alphap = 0$}

In this case constraints \eqref{eq:4e4part_s1} -- \eqref{eq:4e4part_t1} become
\begin{align}
  \label{eq:4e4part_s15}
  &\bvc \cdot \bvcp = 0,
\end{align}
\be
  \label{eq:4_v15}
  \alpha\vec{b} = - \varepsilon_4\vbp\cdot \bvcp - \beta \underline{c} \, \bv{e}_{123} ,
\ee
\begin{align}
  \label{eq:4e4part_v15}
  \alpha \vbp = - \vb \cdot \bvcp - \betap \bvc \, \bv{e}_{123},
\end{align}
\begin{align}
  \label{eq:4e4part_bv15}
  \bvcp = \frac{1}{\alpha} \vbp \wedge \vb,
\end{align}
\be
  \label{eq:4_bv15}
  \alpha \underline{c} = (\varepsilon_4\betap \vbp - \beta \vec{b}) \, \bv{e}_{123},
\ee
\be
  \label{eq:4_t15}
   - \vec{b}\wedge\underline{c} = \alpha \beta \, \bv{e}_{123},
\ee
\begin{align}
  \label{eq:4e4part_t15}
  - \vbp \wedge \bvc = \alpha \betap \, \bv{e}_{123}.
\end{align}
We further have from \eqref{eq:4_v2}, \eqref{eq:4_bv2}, \eqref{eq:4_bv3}, \eqref{eq:4_v3}, \eqref{eq:4e4part_v3} the derived constraints 
\be
  \label{eq:4_v16}
  \vb\cdot \vbp 
  = \beta\betap \, \bv{e}_{123}^2,
\ee
\begin{align}
  \label{eq:4_bv16}
  \bvc^2  
  = - \varepsilon_4 \betap^2 \, \bv{e}_{123}^2 
    + \beta^2 \, \bv{e}_{123}^2,
\end{align}
\begin{align}
  \label{eq:4_v17}
  \vb^2 
  = \varepsilon_4 \bvcp^2 + \beta^2 \, \bv{e}_{123}^2,
\end{align}
\begin{align}
  \label{eq:4e4part_v16}
  \vec{b}^{\prime 2}
  = -\bvcp^2 + \betap^2 \, \bv{e}_{123}^2,
\end{align}
We calculate from \eqref{eq:4_bv15} that
\begin{align}
  \label{eq:4_bv19}
  \alpha^2 \bvc^2 
  &= (\varepsilon_4\betap \vbp - \beta \vb )^2 \, \bv{e}_{123}^2
  \nonumber \\
  &= (\betap^2 \vbpsq + \beta^2 \vb^2 
     - 2 \varepsilon_4 \beta \betap \vbp \cdot \vb)\, \bv{e}_{123}^2
  \nonumber \\
  &\hspace*{-1cm}\stackrel{\eqref{eq:4_v16},\eqref{eq:4_v17},\eqref{eq:4e4part_v16}}{=}   
     [\betap^2 (-\bvcp^2 + \betap^2 \, \bv{e}_{123}^2)
     + \beta^2 (\varepsilon_4 \bvcp^2 + \beta^2 \, \bv{e}_{123}^2)
     - 2 \varepsilon_4 \beta \betap (\beta\betap \, \bv{e}_{123}^2)
     ]\, \bv{e}_{123}^2
  \nonumber \\
  &= \bvcp^2 (-\betap^2 + \varepsilon_4\beta^2) \,\, \bv{e}_{123}^2
     + \betap^4 + \beta^4 -2 \varepsilon_4 \beta^2 \betap^2 
  \nonumber \\
  &=  \bvcp^2 (-\betap^2 + \varepsilon_4\beta^2) \,\, \bv{e}_{123}^2
     + [(-\betap^2 + \varepsilon_4\beta^2)\, \, \bv{e}_{123}^2]^2.
\end{align}
Inserting \eqref{eq:4_bv16} in   \eqref{eq:4_bv19}  we get
\begin{align}
  \label{eq:4_bv18}
  \alpha^2 \bvc^2  
  = \varepsilon_4 \bvc^2 \bvcp^2 + \bvc^4.
\end{align}
If $ \bvc^2 \neq 0$ in \eqref{eq:4_bv18} then
\begin{align}
  \label{eq:4_bv20}
  \varepsilon_4 \alpha^2   
  = \bvcp^2  + \varepsilon_4 \bvc^2,
\end{align}
and the root equation \eqref{eq:4root_a} becomes with \eqref{eq:4_bv16}
\begin{align}
  \label{eq:4root_a18}
4\alpha^2  +3\varepsilon_4 [\bvcp^2 - \varepsilon_4\alpha^2 -\betap^2\, \bv{e}_{123}^2 + \varepsilon_4\beta^2 \, \bv{e}_{123}^2] 
  = 4\alpha^2 = -1,
\end{align}
which has no solution for real $\alpha \neq 0$.

If $\bvc^2 = 0$ the root equation \eqref{eq:4root_a} becomes 
with \eqref{eq:4_bv16} instead
\begin{align}
  \label{eq:4root_a10}
4\alpha^2  +3\varepsilon_4 [\bvcp^2 - \varepsilon_4\alpha^2 ] 
  = \alpha^2 +3\varepsilon_4 \bvcp^2 = -1,
\end{align}
which has again no solution for $\bvcp^2 = 0$.

For $\bvc^2 = 0$ and $\bvcp^2 \neq 0$
\eqref{eq:4e4part_bv1} yields
\begin{align}
  \label{eq:4e4part_bv1_ap0}
  \alpha \bvcp = \vbp \wedge \vb 
  \,\,\, \Longrightarrow \,\,\,
  \alpha^2 \bvcp^2 = (\vbp \wedge \vb )^2 
                   = (\vbp \cdot \vb )^2 - \vbpsq \vb^2 . 
\end{align}
From \eqref{eq:4_v17} and \eqref{eq:4e4part_v16} we can calculate the product $\vbpsq \vb^2$
as
\begin{align}
  \label{eq:b2bp2}
   \vbpsq \vb^2
  &= (-\bvcp^2 + \betap^2 \, \bv{e}_{123}^2)(\varepsilon_4 \bvcp^2 + \beta^2 \, \bv{e}_{123}^2)
  \nonumber \\
  &= -\varepsilon_4 \bvcp^4 
    + \bvcp^2(\underbrace{\varepsilon_4 \betap^2 \, \bv{e}_{123}^2 - \beta^2 \, \bv{e}_{123}^2}_{= -\bvc^2 = 0}) 
    + \beta^2 \betap^2
  \nonumber \\
  &= -\varepsilon_4 \bvcp^4 + \beta^2 \betap^2.  
\end{align}
We now insert \eqref{eq:4_v16} and \eqref{eq:b2bp2} in \eqref{eq:4e4part_bv1_ap0} to obtain
\begin{align}
  \label{eq:4e4part_bv1_ap01}
  \alpha^2 \bvcp^2 
  = \beta^2 \betap^2 +\varepsilon_4 \bvcp^4 
     - \beta^2 \betap^2
  = +\varepsilon_4 \bvcp^4
  \,\,\, \stackrel{\bvcp^2\neq 0}{\Longrightarrow} \,\,\,
  \alpha^2 = \varepsilon_4 \bvcp^2.
\end{align}
Inserting this result into the root equation \eqref{eq:4root_a10} yields again
\be
  4\alpha^2 = -1,
\ee
which as before has no solution for real $\alpha \neq 0$.

\subsection{$n=4$, $\alpha = 0$, $\alphap \neq 0$}

For $\alpha = 0$ the root equation \eqref{eq:4root} simplifies to
\be
  \label{eq:4roota0}
  \vec{b}^2+ \underline{c}^2 + \beta^2\, \bv{e}_{123}^2 
  + \varepsilon_4 \alphap^2 
  - \varepsilon_4 \vbpsq + \varepsilon_4\bvcp^2 - \varepsilon_4\betap^2 \, \bv{e}_{123}^2
  = -1,
\ee
The constraint equations \eqref{eq:4e4part_s1} -- \eqref{eq:4e4part_t1} which have to be satisfied become
\begin{align}
  \label{eq:4e4part_s6}
  &\bvc \cdot \bvcp = 0,
\end{align}
\be
  \label{eq:4_v6}
  \vbp\cdot \bvcp =  - \varepsilon_4 \beta \underline{c} \,\, \bv{e}_{123}  ,
\ee
\begin{align}
  \label{eq:4e4part_v6}
  \vb \cdot \bvcp =  - \betap \bvc \,\, \bv{e}_{123} ,
\end{align}
\be
  \label{eq:4_t6}
   \vec{b}\wedge\underline{c} = 0,
\ee
\begin{align}
  \label{eq:4e4part_t6}
  \vbp \wedge \bvc = 0.
\end{align}
\begin{align}
  \label{eq:4e4part_bv6}
  \alphap \bvc = \vbp \wedge \vb,
\end{align}
\be
  \label{eq:4_bv6}
  \alphap \bvcp = (\betap \vbp - \varepsilon_4\beta \vec{b}) \, \bv{e}_{123} .
\ee
Especially for $\alphap \neq 0$ we obtain from \eqref{eq:4e4part_bv6} and \eqref{eq:4_bv6}
the constraints
\begin{align}
  \label{eq:4e4part_bv7}
  \bvc = \frac{1}{\alphap }\vbp \wedge \vb,
\end{align}
\be
  \label{eq:4_bv7}
  \bvcp = \frac{1}{\alphap }(\betap \vbp - \varepsilon_4\beta \vec{b}) \, \bv{e}_{123} .
\ee
It is obvious that with \eqref{eq:4e4part_bv7} equations \eqref{eq:4_t6} and \eqref{eq:4e4part_t6} are
then fulfilled, because
\be
  \label{eq:bbb0}
  \vb \wedge \vbp \wedge \vb = 0 \text{ and } 
  \vbp \wedge \vbp \wedge \vb = 0 .
\ee
Due to \eqref{eq:bbb0} equation \eqref{eq:4e4part_s6} is also fulfilled
\begin{align}
  \label{eq:4e4part_s7}
  \bvc \cdot \bvcp 
  &\stackrel{\eqref{eq:4e4part_bv7}}{=} \frac{1}{\alphap^2} (\vbp \wedge \vb) \cdot [(\betap \vbp - \varepsilon_4\beta \vb) \, \bv{e}_{123}]
  \nonumber \\
  &= \frac{1}{\alphap^2}[ \betap(\vbp \wedge \vb \wedge \vbp) \, \bv{e}_{123} 
    - \varepsilon_4\beta(\vbp \wedge \vb \wedge \vb)\, \bv{e}_{123}]
  = 0.
\end{align}
Using \eqref{eq:4_bv7} we now check the remaining \eqref{eq:4_v6} and \eqref{eq:4e4part_v6}
\begin{align}
  \label{eq:4_v7}
  \vbp\cdot \bvcp 
  &= \frac{1}{\alphap} \vbp \cdot [(\betap \vbp -  \varepsilon_4 \beta \vb) \, \bv{e}_{123}] 
  \nonumber \\
  &= \frac{1}{\alphap}[\betap \underbrace{\vbp \wedge \vbp}_{=0} \, \bv{e}_{123} 
                      - \varepsilon_4 \beta \underbrace{\vbp \wedge \vb}_{=\alphap \bvc} \, \bv{e}_{123} ] 
  \stackrel{\eqref{eq:4e4part_bv6}}{=} - \varepsilon_4 \beta \bvc \,\, \bv{e}_{123},
  %- \varepsilon_4 \beta \underline{c} \,e_{123}  ,
\end{align}
\begin{align}
  \label{eq:4e4part_v7}
  \vb \cdot \bvcp 
  &= \frac{1}{\alphap}\vb \cdot [(\betap \vbp - \varepsilon_4\beta \vb) \, \bv{e}_{123}]
  \nonumber \\
  &= \frac{1}{\alphap}[\betap \underbrace{\vb \wedge \vbp}_{= - \alphap \bvc} \, \bv{e}_{123} 
                       - \varepsilon_4\beta \underbrace{\vb \wedge \vb}_{=0} \, \bv{e}_{123}]
  \stackrel{\eqref{eq:4e4part_bv6}}{=} - \betap \bvc \, \, \bv{e}_{123} .
  %- \betap \bvc \,\, \bv{e}_{123} ,
\end{align}
Therefore, if the two \textit{constraints} \eqref{eq:4e4part_bv7} 
and \eqref{eq:4_bv7} are satisfied,
all other necessary equations are also satisfied and the 
root equation depends only on
$\alphap$, $\beta$, $\betap$, $\vb$, and~$\vbp:$ 
\begin{align}
  &\vb^2 + \frac{1}{\alphap^2} (\vbp \wedge \vb)^2 
  + \beta^2  \, \bv{e}_{123}^2
  + \varepsilon_4 \alphap^2 
  \nonumber \\
  &- \varepsilon_4 \vb^{\prime 2} 
  + \varepsilon_4 \frac{1}{\alphap^2} (\betap\vbp - \varepsilon_4 \beta \vb )^2 \, \bv{e}_{123}^{2}
  - \varepsilon_4 \betap^2 \, \bv{e}_{123}^{2} = -1.
\end{align}

\subsection{$n=4$, $\alpha = \alphap = 0$}

For $\alpha = \alphap = 0$ the root equation \eqref{eq:4root} simplifies to
\be
  \label{eq:4rootaap0}
  \vec{b}^2+ \underline{c}^2 + \beta^2\, \bv{e}_{123}^2 
  - \varepsilon_4 \vbpsq + \varepsilon_4\bvcp^2 - \varepsilon_4\betap^2 \, \bv{e}_{123}^2
  = -1,
\ee
The constraint equations \eqref{eq:4e4part_s1} -- \eqref{eq:4e4part_t1} 
which have to be satisfied become
\begin{align}
  \label{eq:4e4part_s8}
  &\bvc \cdot \bvcp = 0,
\end{align}
\be
  \label{eq:4_v8}
  \vbp\cdot \bvcp =  - \varepsilon_4 \beta \underline{c} \,\, \bv{e}_{123}  ,
\ee
\begin{align}
  \label{eq:4e4part_v8}
  \vb \cdot \bvcp =  - \betap \bvc \,\, \bv{e}_{123} ,
\end{align}
\be
  \label{eq:4_t8}
   \vec{b}\wedge\underline{c} = 0,
\ee
\begin{align}
  \label{eq:4e4part_t8}
  \vbp \wedge \bvc = 0.
\end{align}
\begin{align}
  \label{eq:4e4part_bv8}
  \vbp \wedge \vb = 0,
\end{align}
\be
  \label{eq:4_bv8}
  \betap \vbp = \varepsilon_4\beta \vec{b}.
\ee

\subsubsection{$n=4$, $\alpha = \alphap = 0$, $\vbp=0$}

For $\vbp=0$ the root equation \eqref{eq:4rootaap0} simplifies to
\be
  \label{eq:4rootaapbp0}
  \vec{b}^2+ \underline{c}^2 + \beta^2\, \bv{e}_{123}^2 
  + \varepsilon_4\bvcp^2 - \varepsilon_4\betap^2 \, \bv{e}_{123}^2
  = -1.
\ee
The remaining constraint equations \eqref{eq:4e4part_s8} -- \eqref{eq:4_bv8} 
which have to be satisfied become
\begin{align}
  \label{eq:4e4part_s9}
  &\bvc \cdot \bvcp = 0,
\end{align}
\be
  \label{eq:4_v9}
  \beta \underline{c} = 0  ,
\ee
\begin{align}
  \label{eq:4e4part_v9}
  \vb \cdot \bvcp =  - \betap \bvc \,\, \bv{e}_{123} ,
\end{align}
\be
  \label{eq:4_t9}
   \vec{b}\wedge\underline{c} = 0,
\ee
\be
  \label{eq:4_bv9}
  \beta \vec{b} = 0.
\ee

\noindent
\textbf{Case: \textbf{$\beta = 0$, $\betap=0$}\\} 
\noindent
Now only the constraints
\begin{align}
  \bvc \cdot \bvcp = 0, 
  \qquad \vb \cdot \bvcp=0, 
  \qquad \vb \wedge \bvc=0
\end{align}
remain.
And the root equation \eqref{eq:4rootaapbp0} reduces to
\begin{align}
  \vb^2 +\bvc^2 + \varepsilon_4 \bvcp^2 = -1.
\end{align}

\noindent
\textbf{Case: \textbf{$\beta = 0$, $\betap \neq 0$}\\} 
\noindent
Now only the constraints
\begin{align}
  \label{eq:aapb0}
  \bvc \cdot \bvcp = 0, 
  \qquad \vb \cdot \bvcp = - \betap \bvc \,\, \bv{e}_{123}, 
  \qquad \vb \wedge \bvc = 0
\end{align}
remain. The second identity in \eqref{eq:aapb0} is equivalent to the \textit{constraint}
\begin{align}
  \label{eq:cond1}
  \bvc = -\frac{1}{\betap}\,\vb\cdot \bvcp \,\, \bv{e}_{123}^{-1}.
\end{align}
We can check that based on \eqref{eq:cond1} the other two constraints of \eqref{eq:aapb0} are also satisfied
\begin{align}
  \label{eq:aapb01}
  \bvc \cdot \bvcp 
  = -\frac{1}{\betap}\,[\vb\cdot \bvcp \,\, \bv{e}_{123}^{-1}]\cdot \bvcp
  = -\frac{1}{\betap}\,[\underbrace{(\vb\cdot \bvcp)\wedge \bvcp}_{=0}]\,\, \bv{e}_{123}^{-1} =0,
  %\qquad \vb \wedge \bvc = 0
\end{align}
and
\begin{align}
  \label{eq:aapb02}
  \vb \wedge \bvc 
  = -\frac{1}{\betap}\,\vb \wedge [\vb\cdot \bvcp \,\, \bv{e}_{123}^{-1}]
  = -\frac{1}{\betap}\,[(\vb \wedge \vb)\cdot \bvcp] \,\, \bv{e}_{123}^{-1} =0.
\end{align}
Inserting $\beta=0$ and \eqref{eq:cond1} into \eqref{eq:4rootaapbp0} yields the root equation
\begin{align}
  \vb^2 + \frac{1}{\betap^2}\,(\vb\cdot \bvcp)^2 \, \bv{e}_{123}^{2} + \varepsilon_4 \bvcp^2 -\varepsilon_4 {\betap}^2\, \bv{e}_{123}^2= -1.
\end{align}

\noindent
\textbf{Case: \textbf{$\beta \neq 0$}\\} 
\noindent
Because of $\beta \neq 0$, the constraints \eqref{eq:4e4part_s9} -- \eqref{eq:4_bv9} reduce to
\begin{align}
  \vb=0, \qquad \bvc = 0.
\end{align}
and the root equation becomes
\begin{align}
  \beta^2 \, \bv{e}_{123}^{2} + \varepsilon_4 \bvcp^2 - \varepsilon_4 {\betap}^2\, \bv{e}_{123}^2 = -1.
\end{align}

\subsubsection{$n=4$, $\alpha = \alphap = 0$, $\vbp\neq0$}
\hspace*{1in}\\

\noindent
\textbf{Case: \textbf{$\vb=0$, $\beta=0$}\\} 
\noindent
This reduces equations \eqref{eq:4e4part_s8} -- \eqref{eq:4_bv8} to the constraints
\begin{align}
  \betap = 0,
  \qquad \bvc \cdot \bvcp = 0, 
  \qquad \vbp\cdot \bvcp=0, 
  \qquad \vbp \wedge \bvc=0,
\end{align}
The root equation \eqref{eq:4rootaap0} becomes then
\begin{align}
  \bvc^2 - \varepsilon_4 {\vb}^{\prime 2} + \varepsilon_4 {\bvcp}^2 = -1.
\end{align}

\noindent
\textbf{Case: \textbf{$\vb=0$, $\beta \neq 0$}\\} 
\noindent
This reduces the constraint equations \eqref{eq:4e4part_s8} -- \eqref{eq:4_bv8} to
\begin{align}
  \label{eq:4e4part_s10}
  &\bvc \cdot \bvcp = 0,
\end{align}
\be
  \label{eq:4_v10}
  \vbp\cdot \bvcp =  - \varepsilon_4 \beta \underline{c} \,\, \bv{e}_{123}  
  \,\,\, \Longrightarrow \,\,\,
  \bvc = -\frac{\varepsilon_4}{\beta} \,\vbp \cdot \bvcp \,\, \bv{e}_{123}^{-1},
\ee
\begin{align}
  \label{eq:4e4part_v10}
  \betap \bvc = 0, 
\end{align}
\begin{align}
  \label{eq:4e4part_t10}
  \vbp \wedge \bvc = 0.
\end{align}
\be
  \label{eq:4_bv10}
  \betap \vbp = 0 
  \,\,\, \Longrightarrow \,\,\,
  \betap  = 0.
\ee
Hence \eqref{eq:4e4part_v10} is satisfied and we must only check \eqref{eq:4e4part_s10}
and \eqref{eq:4e4part_t10}. Inserting \eqref{eq:4_v10} into \eqref{eq:4e4part_s10} gives
\begin{align}
  \label{eq:4e4part_s11}
  \bvc \cdot \bvcp 
  = -\frac{\varepsilon_4}{\beta} \,[\vbp \cdot \bvcp \,\, \bv{e}_{123}^{-1}]\cdot \bvcp
  = -\frac{\varepsilon_4}{\beta} \, [\underbrace{
( \vbp \cdot \bvcp )\wedge \bvcp}_{=0}] \, \bv{e}_{123}^{-1} = 0.
\end{align}
Inserting \eqref{eq:4_v10} into \eqref{eq:4e4part_t10} gives
\begin{align}
  \label{eq:4e4part_t11}
  \vbp \wedge \bvc 
  = -\frac{\varepsilon_4}{\beta} \, \vbp \wedge [\vbp \cdot \bvcp \,\, \bv{e}_{123}^{-1}]
  = -\frac{\varepsilon_4}{\beta} \, [\underbrace{(\vbp \wedge \vbp)}_{=0} \cdot \bvcp] \,\, \bv{e}_{123}^{-1}]
  = 0.
\end{align}
The root equation \eqref{eq:4rootaap0} becomes now with constraints \eqref{eq:4_v10} and \eqref{eq:4_bv10}
\begin{align}
  \frac{1}{\beta^2}(\vbp \cdot \bvcp)^2\, \bv{e}_{123}^{2} + \beta^2 \, \bv{e}_{123}^{2} 
  - \varepsilon_4 {\vb}^{\prime 2} + \varepsilon_4 {\bvcp}^2 = -1.
\end{align}

\noindent
\textbf{Case: \textbf{$\vb\neq0$, $\beta=0$}\\} 
\noindent
This reduces the constraints \eqref{eq:4e4part_s8} -- \eqref{eq:4_bv8} to
\begin{align}
  \label{eq:4e4part_s12}
  &\bvc \cdot \bvcp = 0,
\end{align}
\be
  \label{eq:4_v12}
  \vbp\cdot \bvcp = 0,
\ee
\begin{align}
  \label{eq:4e4part_v12}
  \vb \cdot \bvcp =  - \betap \bvc \,\, \bv{e}_{123} ,
\end{align}
\be
  \label{eq:4_t12}
   \vec{b}\wedge\underline{c} = 0,
\ee
\begin{align}
  \label{eq:4e4part_t12}
  \vbp \wedge \bvc = 0,
\end{align}
\begin{align}
  \label{eq:4e4part_bv12}
  \vbp \wedge \vb = 0
  \,\,\, \Longrightarrow \,\,\,
  \vbp = \gamma \vb, \,\,\,\, 
  \gamma \in \R\setminus \{0\},
\end{align}
\be
  \label{eq:4_bv12}
  \betap \vbp = 0
   \,\,\, \stackrel{\vbp\neq 0}{\Longrightarrow} \,\,\,
   \betap = 0
   \,\,\, \stackrel{\eqref{eq:4e4part_v12}}{\Longrightarrow} \,\,\,
   \vb \cdot \bvcp = 0.
\ee
Constraints \eqref{eq:4e4part_s12} -- \eqref{eq:4_bv12} are equivalent to
($\gamma \in \R\setminus \{0\}$)
\begin{align}
  \bvc \cdot \bvcp = 0, 
  \qquad \vb \wedge \bvc =0,
  \qquad \vb \cdot \bvcp =0,
  \qquad \vbp = \gamma \vb,
  \qquad \betap = 0, 
\end{align}
The root equation \eqref{eq:4rootaap0} then becomes
\begin{align}
  (1-\varepsilon_4\gamma^2)\,\vb^2 + \bvc^2 +\varepsilon_4 \,\bvcp^2 = -1.
\end{align}

\noindent
\textbf{Case: \textbf{$\vb \neq 0$, $\beta \neq 0$}\\} 
\noindent
We obtain from \eqref{eq:4_bv8} that
\be
  \label{eq:4_bv13}
  \betap \vbp = \varepsilon_4\beta \vec{b} 
  \,\,\, \Longrightarrow \,\,\,
  \betap \neq 0 
  \,\,\, \text{ and } \,\,\,
  \vbp = \varepsilon_4\frac{\beta}{\betap}\vec{b},  
\ee
which automatically takes care of \eqref{eq:4e4part_bv8}.
We further calculate from \eqref{eq:4e4part_v8} that
\begin{align}
  \label{eq:4e4part_v13}
  \bvc = -\frac{1}{\betap }\,\vb \cdot \bvcp\,\, \bv{e}_{123}^{-1} .
\end{align}
We now check the remaining four constraints 
\eqref{eq:4e4part_s8}, \eqref{eq:4_v8}, \eqref{eq:4_t8}, \eqref{eq:4e4part_t8} for consistency. 
Due to the proportionality \eqref{eq:4_bv13} of $\vb$ and $\vbp$,  \eqref{eq:4_t8} and \eqref{eq:4e4part_t8} are seen to be equivalent. Inserting \eqref{eq:4e4part_v13} into 
the right hand side of \eqref{eq:4_v8} gives
\begin{align}
  \label{eq:4_v13}
  - \varepsilon_4 \beta \frac{(-1)}{\betap} \vb \cdot \bvcp\,\, \bv{e}_{123}^{-1}\,\, \bv{e}_{123}  
  = \varepsilon_4 \frac{\beta}{\betap} \vb \cdot \bvcp
  \stackrel{\eqref{eq:4_bv13}}{=} \vbp\cdot \bvcp.
\end{align}
Inserting \eqref{eq:4e4part_v13} into \eqref{eq:4e4part_s8} gives
\begin{align}
  \label{eq:4e4part_s13}
  \bvc \cdot \bvcp
  = -\frac{1}{\betap }\,[\vb \cdot \bvcp\,\, \bv{e}_{123}^{-1}] \cdot \bvcp
  = -\frac{1}{\betap }\,[\underbrace{(\vb \cdot \bvcp)\wedge \bvcp}_{=0}]\,\, \bv{e}_{123}^{-1}
  =0.
\end{align}
Finally inserting \eqref{eq:4e4part_v13} into \eqref{eq:4_t8} gives
\begin{align}
  \label{eq:4_t13}
   \vb \wedge \bvc 
   = -\frac{1}{\betap } \,\vb \wedge [\vb \cdot \bvcp\,\, \bv{e}_{123}^{-1}]
   = -\frac{1}{\betap } \,[\underbrace{(\vb \wedge \vb)}_{=0} \cdot \bvcp]\,\, \bv{e}_{123}^{-1}
   =0 .
\end{align}
Everything is therefore consistent and with the constraints \eqref{eq:4_bv13} and 
\eqref{eq:4e4part_v13} for $\vbp$ and $\bvc$ we get from \eqref{eq:4rootaap0} the root equation
\begin{align}
  (1 - \varepsilon_4 \frac{\beta^2}{\betap^2})\vb^2 + \frac{1}{\betap^2}(\vb \cdot \bvcp)^2\, \bv{e}_{123}^{2}
  + \beta^2 \, \bv{e}_{123}^{2}
  + \varepsilon_4 \bvcp^2 
  - \varepsilon_4 \betap^2 \, \bv{e}_{123}^{2} = -1.
\end{align}
This concludes the discussion of $n=4$, $\alpha = \alphap =0$.

Table \ref{tb:roots4} on page  \pageref{tb:roots4} lists all geometric
roots of $-1$ of $\G_{p,q}$, $n=p+q=4$. We point out, that similar to Table 
\ref{tb:roots123}, also in Table \ref{tb:roots4} all root equations of the third column result from the general $n=4$ root equation \eqref{eq:4root}, simply by inserting the case conditions and constraints of columns one and two.

%\section{Relation with idempotents and group theory}

\section{Conclusions}

\label{pg:roots123}
\begin{table}%[ht]
\begin{center}
\caption{Geometric roots of $-1$ for Clifford algebras
$\G_{p,q}$, $n=p+q\leq 3$. The multivectors are denoted for $n=1$ by
$\alpha + \beta \vect{e}_1$, 
for $n=2$ by
$\alpha + b_1 \vect{e}_1 + b_2 \vect{e}_2 + \beta \bv{e}_{12}$, 
and for $n=3$ by
$\alpha 
  + b_1 \vect{e}_1 + b_2 \vect{e}_2 + b_3 \vect{e}_3
  + c_{1} \bv{e}_{23} + c_{2} \bv{e}_{31} + c_{3} \bv{e}_{12} 
  + \beta \bv{e}_{123} $. %\vspace*{2mm}
\label{tb:roots123}}
\begin{tabular}{lll}
%\hline\noalign{\smallskip}
\topline
%Crystal&\multicolumn{2}{c}{Oblique}
$n$ & Cases & Solutions $A$ and root equations
\\
\midline
$1$& & no solution for $\G_1$  \\
& & $A=\pm \bv{e}_1$ for $\G_{0,1}$ 
\\
\midline
\rule{0mm}{4.0mm}%
$2$ & $\alpha = 0$ & 
$\beta^2 
= b_1^2 \varepsilon_2 + b_2^2 \varepsilon_1 +\varepsilon_1\varepsilon_2$
\\
&&
$\hspace*{0.0mm} \beta^2= \left\{ 
  \begin{array}{cl}
  \phantom{-}b_1^2 + b_2^2 +1  & \mbox{ for } \G_2\\
  -b_1^2 + b_2^2 -1 & \mbox{ for } \G_{1,1}\\
  -b_1^2 - b_2^2 +1 & \mbox{ for } \G_{0,2}
  \end{array} 
\right. 
$\rule[-7mm]{0mm}{15mm}% 
\\
\cline{2-3}%\noalign{\smallskip}
\rule{0mm}{5.0mm}%
 & $\alpha \neq 0$ & no solution \\
\midline
$3$ & Constraint: & $\phantom{-}0=\vect{b}\wedge\bvc = b_1c_1+b_2c_2+b_3c_3$ \\
$$ & $\alpha = \beta = 0$ & $-1 = \vec{b}^2+ \underline{c}^2 $ \\
& $\quad$  & $\hspace*{0.0mm}-1= b_1^2\varepsilon_1 +  b_2^2\varepsilon_2 + b_3^2\varepsilon_3
- c_1^2 \varepsilon_2\varepsilon_3 - c_2^2 \varepsilon_3\varepsilon_1 - c_3^2 \varepsilon_1\varepsilon_2$ \\
%& $\qquad \vect{b}\wedge\bvc=0$ & $$
%& $\qquad \vect{b}\wedge\bvc=0$ & $\hspace*{11mm}- c_1^2 \varepsilon_2\varepsilon_3 - c_2^2 \varepsilon_3\varepsilon_1 - c_3^2 \varepsilon_1\varepsilon_2$
%
%\\
&&
  $
  \hspace*{0.0mm}-1= \left\{
    \begin{array}{ll}
     \phantom{-(} b_1^2 +  b_2^2 + b_3^2\phantom{)} - (c_1^2 + c_2^2 + c_3^2) & \text{for } \G_{3}\\
     \phantom{-(} b_1^2 -  b_2^2 - b_3^2\phantom{)} - (c_1^2 - c_2^2 - c_3^2) & \text{for } \G_{1,2} \\
     \phantom{-(} b_1^2 +  b_2^2 - b_3^2\phantom{)} + (c_1^2 + c_2^2 - c_3^2) & \text{for } \G_{2,1}\\
      - (b_1^2 +  b_2^2 + b_3^2) - (c_1^2 + c_2^2 + c_3^2)\hspace*{-2mm} & \text{for } \G_{0,3}
    \end{array}
  \right. $\rule[-9mm]{0mm}{20mm}%
\\
\cline{2-3}%\noalign{\smallskip}
\rule{0mm}{4.0mm}%
& $\alpha = 0$, $\beta \neq 0$ 
& $A = \pm \bv{e}_{123}$ for $\G_3$, $\G_{1,2}$  \\
&%\multicolumn{1}{c}{a}
& 
no solution for $\G_{2,1}$, $\G_{0,3}$\rule[-2mm]{0mm}{6mm}% 
\\
\cline{2-3}%\noalign{\smallskip}
\rule{0mm}{4.0mm}%
&$\alpha \neq 0$& no solution \\
\bottomline
\end{tabular}
\end{center}
\end{table}

\begin{table}%[ht]
\begin{center}
\caption{Geometric roots of $-1$ for Clifford algebras
$\G_{p,q}$, $n=p+q = 4$. The multivectors are denoted by
$\alpha + \vb + \bvc + \beta \, \bv{e}_{123} 
   + (\alphap + \vbp + \bvcp + \betap \, \bv{e}_{123})\,\vect{e}_4$,
for details see \eqref{eq:4multiv} in the text.
\label{tb:roots4}}
\begin{tabular}{l|l|l}
%\hline\noalign{\smallskip}
\hline
\rule{0mm}{5.0mm}% 
%Crystal&\multicolumn{2}{c}{Oblique}
Case & Subcase / Constraints & Solutions and root equations
\\
\hline
\rule{0mm}{5.0mm}% 
$\alpha \neq 0$ &  & no solution \\
\hline %\noalign{\smallskip}
\rule{0mm}{5.0mm}%
$\alpha = 0$, 
  & $\quad \bvc = \frac{1}{\alphap }\vbp \wedge \vb,$ 
  & $\vb^2 + \frac{1}{\alphap^2} (\vbp \wedge \vb)^2 
  + \beta^2  \, \bv{e}_{123}^2
  + \varepsilon_4 \alphap^2 
  - \varepsilon_4 \vb^{\prime 2}$ \\
$\alphap \neq 0$
  & $\quad \bvcp \hspace*{-0.8mm}= \hspace*{-0.8mm}\frac{1}{\alphap }(\betap \vbp \hspace*{-1mm}- \hspace*{-1mm}\varepsilon_4\beta \vec{b}) \, \bv{e}_{123}\hspace*{-1.5mm}$ 
  & $\,\,\, + \varepsilon_4 \frac{1}{\alphap^2} (\betap\vbp - \varepsilon_4 \beta \vb )^2 \, \bv{e}_{123}^{2}
  - \varepsilon_4 \betap^2 \, \bv{e}_{123}^{2}$ \\
&& $\,\,\, = -1$ \\
\hline %\noalign{\smallskip}
\rule{0mm}{5.0mm}%
$\alpha = 0$,  & $\beta = \betap =0$ & \\
$\alphap = 0$, & $\quad \bvc \cdot \bvcp = 0,$ 
& $\vb^2 +\bvc^2 + \varepsilon_4 \bvcp^2 = -1$ \\
$\vbp=0$       & $\quad
   \vb \cdot \bvcp=0, \,\,
   \vb \wedge \bvc=0$ 
& \\
\cline{2-3}%\noalign{\smallskip}
\rule{0mm}{5.0mm}%
& $\beta = 0,$ $\betap \neq 0$ & \\
& $\quad \bvc = -\frac{1}{\betap}\,\vb\cdot \bvcp \,\, \bv{e}_{123}^{-1}$ 
& $ \vb^2 + \frac{1}{\betap^2}\,(\vb\cdot \bvcp)^2 \, \bv{e}_{123}^{2} + \varepsilon_4 \bvcp^2 -\varepsilon_4 {\betap}^2\, \bv{e}_{123}^2$ \\
&& $ \,\,\, = -1$ \\
\cline{2-3}%\noalign{\smallskip}
\rule{0mm}{5.0mm}%
& $\beta \neq 0$ & \\
& $\quad \vb=0,$ \,\, $\bvc=0$ 
& $\beta^2 \, \bv{e}_{123}^{2} + \varepsilon_4 \bvcp^2 - \varepsilon_4 {\betap}^2\, \bv{e}_{123}^2 = -1$ \rule[-3mm]{0mm}{8mm}%
\\
\hline %\noalign{\smallskip}
\rule{0mm}{5.0mm}%
$\alpha = 0,$  & $\vb=0,$ $\beta =0$ & \\
$\alphap = 0,$ & $\quad \betap = 0$, \,\, $\bvc \cdot \bvcp = 0,$  
& $\bvc^2 - \varepsilon_4 {\vb}^{\prime 2} + \varepsilon_4 {\bvcp}^2 = -1$\\
$\vbp \neq 0$ & $\quad \vbp \cdot \bvcp = 0$, \,\, $\vbp \wedge \bvc = 0\hspace*{-2mm}$ & \\
\cline{2-3}%\noalign{\smallskip}
\rule{0mm}{5.0mm}%
& $\vb=0,$ $\beta \neq 0$ & \\
& $\quad  \bvc = -\frac{\varepsilon_4}{\beta} \,\vbp \cdot \bvcp \,\, \bv{e}_{123}^{-1},$  
& $\frac{1}{\beta^2}(\vbp \cdot \bvcp)^2\, \bv{e}_{123}^{2} + \beta^2 \, \bv{e}_{123}^{2} 
  - \varepsilon_4 {\vb}^{\prime 2} + \varepsilon_4 {\bvcp}^2$\\
& $\quad \betap = 0$ & $ \,\,\, = -1 $ \\
\cline{2-3}%\noalign{\smallskip}
\rule{0mm}{5.0mm}%
& $\vb \neq 0,$ $\beta =0$ & \\
&  $\quad \bvc \cdot \bvcp = 0,$ \,\, $\betap = 0, $ & \\
&  $\quad \vb \wedge \bvc =0,$ \,\, $\vb \cdot \bvcp =0,$ 
&  $(1-\varepsilon_4\gamma^2)\,\vb^2 + \bvc^2 +\varepsilon_4 \,\bvcp^2 = -1$ \\
\rule[-2mm]{0mm}{5.9mm}% 
&  $\quad \vbp = \gamma \vb, \,\,\, \gamma \in \R\setminus \{0\}$ & \\
\cline{2-3}%\noalign{\smallskip}
\rule{0mm}{5.0mm}%
& $\vb \neq 0,$ $\beta \neq 0$ & \\
& $\quad \betap \neq 0, \,\,\, \vbp = \varepsilon_4\frac{\beta}{\betap}\vec{b},$ 
& $(1 - \varepsilon_4 \frac{\beta^2}{\betap^2})\vb^2 + \frac{1}{\betap^2}(\vb \cdot \bvcp)^2\, \bv{e}_{123}^{2} + \beta^2 \, \bv{e}_{123}^{2}
  $ \\
\rule[-3mm]{0mm}{5.9mm}% 
& $\quad \bvc = -\frac{1}{\betap }\,\vb \cdot \bvcp\,\, \bv{e}_{123}^{-1}$ 
& $ \quad  + \varepsilon_4 \bvcp^2 
  - \varepsilon_4 \betap^2 \, \bv{e}_{123}^{2} = -1$ 
\\
\hline
\end{tabular}
\end{center}
\end{table}

Table \ref{tb:roots123} lists all geometric roots of $-1$ for Clifford algebras $\G_{p,q}$, 
$n=p+q\leq 3$, and Table \ref{tb:roots4} does the same for Clifford algebras $\G_{p,q}$, 
$n=p+q = 4$. The content of both tables has been checked with the MAPLE package 
CLIFFORD~\cite{AF:CLIFFORD}. The solutions for $\G_3$ included in Table~\ref{tb:roots123} 
correspond to the biquaternion roots of $-1$ found in~\cite{SJS:Biqroots}.

Overall the calculations and the results demonstrate how in Clifford algebras extensive calculations can be done without referring to coordinates \cite{HS:CAtoGC,AOR}. In the case of $\G_{p,q}$, 
$n=p+q = 4,$ we arbitrarily selected one non-isotropic vector $\vec{e}_4$ for suitably splitting the algebra in order to use well developed techniques for algebras $\G_{p,q}$, $n=p+q = 3$. 
In the end it is always possible to express the results in coordinates as in 
Table~\ref{tb:roots123}. However, this considerably blows up the expressions
and blurs the mostly $p,q$-signature independent form of the root equations 
of the families of geometric roots of $-1$. In the case of $\G_{p,q}$, $n=p+q = 4,$ in 
Table~\ref{tb:roots4}, we have not expressed the results in coordinates, because then the table would extend over several pages. 

Open questions are:
\begin{itemize}
\item
How can the graded structure of $\G_{p,q}$ be used best in the calculation of higher order geometric multivector square roots of $-1$? This also includes a question how to best use, for this type of computation, invariance of the equation $AA =-1$ under Clifford algebra (anti) automorphisms such as grade involution, reversion or conjugation, and under symmetries of the root equation. For example, under the grade involution,
\be
  AA = -1
  \,\,\, \Longleftrightarrow \,\,\,
  \hat{A}\hat{A} = -1.
\ee
Another example would be a rotor $R$ symmetry
\begin{align}
  & AA = \langle AA \rangle = -1
  \,\,\, \Longleftrightarrow \,\,\, \\
  & R^{-1}ARR^{-1}AR = \langle R^{-1}ARR^{-1}AR \rangle 
  = \langle R^{-1}AAR \rangle =\langle AA \rangle= -1.
  \nonumber 
\end{align}
\item
The interesting relationship with families of idempotents of Clifford geometric
algebras~\cite{AFPR:idem}.
\item
What is the relationship with combinatorics?
\item
Expansion of this work to Clifford algebras $\G_{p,q}$ in arbitrary dimension $n=p+q$. 
For this purpose, it will be appropriate to use the modulo eight periodicity of Clifford algebras and the isomorphisms with matrix rings. Central elements squaring to $-1$ would be of particular importance as then they can be used in place of the imaginary $i.$
\item
The further use of Clifford algebra computation software like CLIFFORD for MAPLE and other packages~\cite{AS:ICCA6lec,AF:CLIFFORD,EH:CliffSoft}.  
\end{itemize}

Of special interest in physics are the Clifford algebras of Minkowski space-time,
sometimes called~\cite{DH:STA} \textit{space-time algebras} $\G_{3,1}$ and $\G_{1,3}$. 
Table~\ref{tb:roots4} contains the complete set of all geometric roots of $-1$ for these algebras, so in particular all possible geometric multivector elements that may take on the role of the imaginary unit $i$ in quantum mechanics, which is e.g. fundamental for the description of spin and for wave propagation.

Finally, the door is now wide open to construct all possible new types of 
Clifford Fourier transformations (CFT)~\cite{SBS:FastCmQFT} for multivector fields with domains and image domains ranging over the full Clifford algebras involved or subalgebras and subspaces thereof. In particular all known Fourier transformations will find their place in this new
general framework. The close relationship of wavelet transformations~\cite{MH:CliffWavl} and 
windowed transformations~\cite{MHA:2DCliffWinFT} to Fourier transformations shows that also in these fields new mathematics is to be expected. 

Examples of CFTs working with non-central replacements of the imaginary unit $i$ are the quaternion FT (QFT)~\cite{GS:ncHCFT,TB:thesis,EH:QFTgen,MHAV:WQFT}, and the CFT~\cite{ES:CFTonVF,HM:CFToMVF} where $i$ is replaced by pseudoscalars in $\G_{n}, n=2 \, (\rm mod \, 4)$. This shows that in principle every geometric root of $-1$, be it central or not, gives rise to its own geometric FT. Regarding the non-central geometric roots of $-1$, the example of the QFT shows that the non-commutativity may indeed be of advantage for obtaining more information about the symmetry and the physical nature of signals thus processed.

\section*{Acknowledgments}

E.H. thanks God, the Creator:
\textit{How great are your works, O LORD,
       how profound your thoughts!} \cite{NIV:Psalm92v5}.
He thanks his family for their total loving support.
We gratefully acknowledge valuable advice given by H. Ishii (Nagoya) and by
J. Helmstetter (Grenoble).


\begin{thebibliography}{99}
\bibitem{SJS:Biqroots}
S. J. Sangwine,
\textit{Biquaternion (Complexified Quaternion)
Roots of -1},
Adv. Appl. Cliford Alg. \textbf{16}(1), pp. 63–-68, 2006.

\bibitem{WKC:AppGrass}
W.K. Clifford
\textit{Applications of Grassmann's Extensive Algebra},
American Journal of Mathematics Pure and Applied \textbf{1}, pp. 350--358, 1878; 
Online: http://sinai.mech.fukui-u.ac.jp/gcj/wkconline.html


\bibitem{HS:CAtoGC}
D. Hestenes, G. Sobczyk,
\textit{Clifford Algebra to Geometric Calculus}, 
Kluwer, 1984.

          
\bibitem{BDS:CA} 
F. Brackx, R. Delanghe, and F. Sommen, 
\textit{Clifford Analysis}, 
Vol. 76 
of Research Notes in Mathematics, Pitman Advanced Publishing Program, 1982.


\bibitem{GS:ncHCFT} 
T. B\"{u}low, M. Felsberg and G. Sommer,
\textit{Non-commutative Hypercomplex Fourier Transforms 
of Multidimensional Signals}, 
in G. Sommer (ed.),
\textit{Geom. Comp. with Cliff. Alg.,
Theor. Found. and Appl. in Comp. Vision and Robotics,} 
Springer (2001), 187--207.


\bibitem{LMQ:CAFT94} 
C. Li, A. McIntosh and T. Qian, 
\textit{Clifford Algebras, Fourier Transform
and Singular Convolution Operators On Lipschitz Surfaces}, 
Revista Matematica Iberoamericana,
\textbf{10 (3)}, (1994), 665--695.


\bibitem{AM:CAFT96}
A. McIntosh,
\textit{Clifford Algebras, Fourier Theory, Singular Integrals, and 
Harmonic Functions on Lipschitz Domains}, 
chapter 1 of J. Ryan (ed.), 
Clifford Algebras in Analysis and Related Topics, 
CRC Press, Boca Raton, 1996.


\bibitem{TQ:PWT}
T. Qian,
\textit{Paley-Wiener Theorems and Shannon Sampling
in the Clifford Analysis Setting}
in R.~Ab\l amowicz (ed.), Clifford Algebras - 
Applications to Mathematics, Physics, and Engineering,
Birk\"{a}user, Basel, (2004), 115--124.


\bibitem{ES:CFTonVF} 
J. Ebling and G. Scheuermann,
\textit{Clifford Fourier Transform on Vector Fields}, 
IEEE Transactions on Visualization and Computer Graphics, \textbf{11 (4)}, 
July/August (2005), 469--479.


\bibitem{HM:CFTUP}
B. Mawardi, E. Hitzer, 
\textit{Clifford Fourier Transformation and
Uncertainty Principle for the Clifford algebra $Cl_{3,0}$},
Adv. Appl. Clifford Alg., \textbf{16}(1) (2006), pp. 41--61.


\bibitem{EM:CFaUP} 
E. Hitzer, B. Mawardi,
\emph{Uncertainty Principle 
for the Clifford algebra $ Cl_{n,0}$, n = 3(mod 4)
based on Clifford Fourier transform,}  
in Springer SCI book series \textit{Applied and
Numerical Harmonic Analysis}, 2006, pp. 45--54.


\bibitem{HM:CFToMVF} 
E. Hitzer, B. Mawardi,
\textit{Clifford Fourier Transform on Multivector Fields and 
Uncertainty Principles
for Dimensions $n=2 \, (\rm mod \, 4)$ and $n=3 \, (\rm mod \, 4)$},
Adv. Appl. Clifford Alg., \textbf{18}(3-4), (2008), pp. 715--736.

\bibitem{EH:QFTgen}
E. Hitzer,
\textit{Quaternion Fourier Transform on Quaternion
Fields and Generalizations}
Adv. Appl. Clifford Alg. \textbf{17}(3), (2007), pp. 497--517.


\bibitem{TB:thesis}
T. B\"{u}low, 
\textit{Hypercomplex Spectral Signal Representations for 
the Processing and Analysis of Images},
PhD Thesis, Univ. of Kiel, 1999. 


\bibitem{MF:thesis} M. Felsberg, 
\textit{Low-Level Image Processing with the Structure
Multivector}, PhD thesis, Univ. of Kiel, 2002.


\bibitem{DH:NF1}
D. Hestenes,
\textit{New Foundations for Classical Mechanics}, 
Kluwer, 1999.


\bibitem{DH:STA}
D. Hestenes,
\textit{Space-Time Algebra},
Gordon and Breach, 1966.

\bibitem{HL:IAGR}
H. Li, 
\textit{Invariant Algebras and Geometric Reasoning},
World Scientific, Singapore, 2009.

\bibitem{LD:InP}
L. Dorst, 
\textit{The Inner Products of Geometric Algebra},
in L. Dorst, C. Doran, J. Lasenby (eds.), Birkh\"{a}user,
Boston, pp. 35-46, 2002.


\bibitem{PL:CAaS}
P. Lounesto, 
\textit{Clifford Algebra and Spinors}, 
2nd ed., CUP,
Cambridge, 2006.


\bibitem{AFPR:idem}
R.~Ab\l amowicz, B. Fauser, K. Podlaski, J. Rembieli\'{n}ski,
\textit{Idempotents of Clifford Algebras},
Czechoslovak Journal of Physics, \textbf{53}(11), pp. 949--954, 2003.


\bibitem{AF:CLIFFORD}
R.~Ab\l amowicz, B. Fauser,
\textit{CLIFFORD - A Maple Package for Clifford Algebra Computations with Bigebra, SchurFkt, GfG - Groebner for Grassmann, Cliplus, Define, GTP, Octonion, SP, SymGroupAlgebra, and code\_support}, http://www.math.tntech.edu/rafal/, December 2008.

\bibitem{AOR} 
R.~Ab\l amowicz, Z.~Oziewicz and J.~Rzewuski, 
\textit{Clifford algebra approach to twistors}, 
J. Math. Phys. \textbf{23} (1982), 231--242.


\bibitem{AS:ICCA6lec}
R.~Ab\l amowicz, G. Sobczyk,
Appendix 7.1 of
\textit{Lectures on Clifford (Geometric) Algebras and Applications},
Birkh\"{a}user, Boston, 2004.


\bibitem{EH:CliffSoft}
E. Hitzer,
Geometric Calculus International -- Software, 
http://sinai.mech.fukui-u.ac.jp/gcj/gc\_int.html\#software

\bibitem{SBS:FastCmQFT}
S. Said, N. Le Bihan, S. J. Sangwine,
\textit{Fast complexified quaternion Fourier transform},
IEEE Transactions on Signal Processing, \textbf{56}(4), pp. 1522--1531, 2008.

\bibitem{MH:CliffWavl}
B. Mawardi, E. Hitzer,
\textit{Clifford Algebra $Cl(3,0)$-valued Wavelet Transformation, Clifford Wavelet Uncertainty Inequality and Clifford Gabor Wavelets}, 
International Journal of Wavelets, Multiresolution and Information Processing, 
\textbf{5}(6), pp. 997-1019, 2007.

\bibitem{MHA:2DCliffWinFT}
B. Mawardi, E. Hitzer, S. Adji, 
\textit{Two-Dimensional Clifford Windowed Fourier Transform}, 
accepted for G. Scheuermann, E. Bayro-Corrochano (eds.), Applied Geometric Algebras in
Computer Science and Engineering, Springer, New York, 2009.

\bibitem{MHAV:WQFT}
B. Mawardi, E. Hitzer, R. Ashino, R. Vaillancourt, 
\textit{Windowed Fourier transform of two-dimensional quaternionic signals}, submitted to Appl. Math. and Computation, March 2009.


\bibitem{NIV:Psalm92v5}
Psalm 92, verse 5, 
\textit{New Int. Version of the Bible}, 
www.biblegateway.com

\end{thebibliography}
\end{document}